\documentclass
[11pt, a4paper]{article}

\usepackage{inputenc} 
\usepackage[english]{babel}
\usepackage{dsfont}
\usepackage{authblk}

\usepackage[totalheight=24 true cm, totalwidth=15 true cm]{geometry}
\usepackage{enumitem}

\setlist[itemize]{noitemsep}
\setlist[enumerate]{noitemsep}
\usepackage{bbding}
\usepackage{amsmath,amsthm,amssymb}
\usepackage{mathrsfs}
\usepackage{wasysym}
\usepackage{mathtools}
\usepackage{stmaryrd}
\usepackage{url}
\usepackage{wrapfig}
\usepackage{starfont}
\usepackage{pifont}
\usepackage{eurosym}
\usepackage{subcaption}
\usepackage{setspace}
\usepackage{bbm}
\usepackage{dsfont}
\usepackage{url}
\usepackage{amsmath,blkarray,booktabs, bigstrut}
\usepackage{multirow,multicol}
\usepackage{xcolor}
\usepackage{color}
\usepackage{imakeidx}
\makeindex[]

\usepackage{algorithm}
\usepackage{algpseudocode}
\algrenewcommand\algorithmicrequire{\textbf{Input:}}
\algrenewcommand\algorithmicensure{\textbf{Output:}}
\usepackage{relsize}
\usepackage{float}
\usepackage{tikz,pgf}
\usepackage{tikz-cd}
\usepackage{wasysym}
\usepackage{graphicx}
\usepackage{fancyhdr}
\usepackage{stmaryrd}
\usepackage{verbatim}
\usepackage{pgfplots}
\pgfplotsset{compat=1.15}
\usetikzlibrary{arrows}
\usepackage{xspace}

\usepackage[all,dvips,knot,web,arc,curve,color,frame]{xy}
\usepackage[all,knot,arc,color,web]{xy}
\xyoption{arc}
\xyoption{web}
\xyoption{curve}

\numberwithin{equation}{section}

\theoremstyle{plain}
\newtheorem{theorem}{Theorem}[section]
\newtheorem{proposition}[theorem]{Proposition}
\newtheorem{corollary}[theorem]{Corollary}
\newtheorem{lemma}[theorem]{Lemma}
\newtheorem{notation}[theorem]{Notation}

\theoremstyle{definition}
\newtheorem{definition}[theorem]{Definition}

\newtheorem{remark}[theorem]{Remark}

\newtheorem{example}[theorem]{Example}
\newtheorem{question}{Question}

\newtheorem{conjecture}[theorem]{Conjecture}

\usepackage[bookmarksnumbered=true]{hyperref} 
\hypersetup{
     colorlinks = true,
     linkcolor = blue,
     anchorcolor = blue,
     citecolor = teal,
     filecolor = blue,
     urlcolor = blue
     }
\newcommand\restr[2]{{
  \left.\kern-\nulldelimiterspace 
  #1 
  \vphantom{\big|} 
  \right|_{#2} 
  }}

\newcommand{\pms}{\mathcal{P}_{\mS'}}

\everymath{\displaystyle}
\allowdisplaybreaks

\usepackage{graphicx}

\newcommand{\size}[1]{\left| #1 \right|} 

\newcommand{\ini}{\textup{in}_{<}}

\newcommand{\VPJ}{\textup{v}_{P_{\mathcal{S}}}(J_G)}
\newcommand{\VPT}{\textup{v}_{P_{\mathcal{S}}}(J_T)}

\newcommand{\DCS}{\mathcal{D}_{c}^{-}(\mathcal{S})}
\newcommand{\ECS}{\mathcal{D}_{c}^{+}(\mathcal{S})}
\newcommand{\JS}{J_{\mathcal{T}(\mathcal{S})}}
\newcommand{\pJS}{J_{\mathcal{T}'(\mathcal{S})}}

\newcommand{\TS}{\mathcal{T}(\mathcal{S})}
\newcommand{\pTS}{\mathcal{T}'(\mathcal{S})}

\newcommand{\reg}{\textup{reg}}

\newcommand{\vn}{\mathrm{v}}

\newcommand{\lvs}{\textup{v}_{\mathcal{S}}^{-}(G)}
\newcommand{\lvt}{\textup{v}_{\mathcal{S}}^{-}(T)}

\newcommand{\pti}{\mathcal{T}'(\{i\})}
\newcommand{\uvs}{\textup{v}_{\mathcal{S}}^{+}(G)}
\newcommand{\uvt}{\textup{v}_{\mathcal{S}}^{+}(T)}

\newcommand{\iv}{\textup{iv}}

\DeclareMathOperator{\supp}{supp}

\DeclareMathOperator{\mS}{\mathcal{S}}
\DeclareMathOperator{\CC}{\mathbb{C}}

\usepackage[maxbibnames=99]{biblatex}
\hypersetup{
  pdftitle={Combinatorial bounds for the v-number of binomial edge ideals},
  pdfauthor={Emiliano Liwski}
}

\newtheorem{theoremA}{Theorem}

\newtheorem{theoremB}{Theorem}

\newtheorem{theoremC}{Theorem}

\title{\vspace*{-2em}\textsc{Combinatorial bounds for the v-number of binomial edge ideals}}

\author{Emiliano Liwski}

\date{\vspace*{-.5em}\today\vspace*{-1.5em}}

\begin{document}
\maketitle

\begin{abstract}
We study the local $\mathrm v$-numbers of binomial edge ideals
from a combinatorial perspective.
For each minimal prime $P_{\mathcal S}$ of $J_G$, indexed
by a cutset $\mathcal S$ of $G$, we consider two combinatorial
invariants, $\mathrm v^-_{\mathcal S}(G)$ and
$\mathrm v^+_{\mathcal S}(G)$, defined in terms of transversals
of suitable hypergraphs, and show that the local $\mathrm v$-number
$\mathrm v_{P_{\mathcal S}}(J_G)$ satisfies
$\mathrm v^-_{\mathcal S}(G)
\leq \mathrm v_{P_{\mathcal S}}(J_G)
\leq \mathrm v^+_{\mathcal S}(G)$.
Under a suitable hypothesis on the connected components,
we provide explicit graph-theoretic descriptions of these
invariants and analyze when they coincide, giving a
combinatorial way to determine the local $\mathrm v$-numbers.
As an application, we derive a sufficient condition for this coincidence, which applies to a broad class of graphs. Furthermore, we identify a constant that controls the gap between $\lvs$ and $\uvs$, thereby showing that the local $\mathrm{v}$-number lies within a bounded interval. Applications to cycles and trees illustrate both the
equality criterion and the resulting bounds.

\end{abstract}

\textbf{MSC:} 05E40, 13P10, 13P25, 13C05, 05C69, 94B27

\textbf{Keywords:} binomial edge ideals; v-number; graph invariants; connected domination; cutsets; transversals.

\section{Introduction}

Binomial edge ideals were introduced independently by Herzog, Hibi, Hreinsdóttir, Kahle, and Rauh \cite{herzog2010binomial} and by Ohtani \cite{ohtani2011graphs}. They connect
graph theory with commutative algebra and arise in algebraic
statistics through conditional independence statements \cite{ambhore2024v,dey2024v,herzog2010binomial,ohtani2011graphs}. 


\begin{definition}\label{definition binomial edge}
Let $S=\mathbb{K}[x_1,\ldots, x_n, y_1,\ldots,y_n]$, and let $G$ be a simple graph with vertex set $V(G)=[n]$ and edge set $E(G)$. The {\em binomial edge ideal} of $G$ is the ideal

$$
J_G:=(f_{i,j}:=x_i y_j-x_j y_i:\{i,j\}\in E(G),\ i<j).
$$

\end{definition}

For the results of this paper, we take $\mathbb K=\mathbb C$. The combinatorics of $G$ is reflected in several algebraic invariants of $J_G$. In this paper, we focus on the $\mathrm{v}$-number. We write $\operatorname{Ass}(I)=\operatorname{Ass}(S/I)$ and denote
the degree-$d$ homogeneous component of $S$ by $S_d$.

\begin{definition}\textup{\cite[Definition~4.1]{cooper2020generalized}}\label{first definition}
Given a homogeneous ideal $I$ of $S=\mathbb{K}[x_1,\ldots,x_m]$, the $\mathrm{v}$-number of $I$ is

$$
\mathrm{v}(I):=\min\{d\geq 0:\exists f\in S_d\text{ and }\mathfrak{p}\in\operatorname{Ass}(I)
\text{ with }(I:f)=\mathfrak{p}\}.
$$

For each associated prime $\mathfrak{p}$ of $I$, the local $\mathrm{v}$-number is defined by

$$
\mathrm{v}_{\mathfrak{p}}(I):=\min\{d\geq 0:\exists f\in S_d\text{ with }(I:f)=\mathfrak{p}\}.
$$

Then
$\mathrm{v}(I)=\min\{\mathrm{v}_{\mathfrak{p}}(I):\mathfrak{p}\in\operatorname{Ass}(I)\}.
$
\end{definition}

The v-number has been studied through its combinatorial and
homological properties; see \cite{saha2024cover,biswas2024study,civan2023v,saha2024v,saha2022v,jaramillo2021v}. For squarefree
monomial ideals, \cite{jaramillo2021v} gives a combinatorial formula with an
interpretation in terms of an associated simplicial complex.
Several works also compare the v-number with
Castelnuovo--Mumford regularity. In particular, for certain classes
of graphs, the monomial edge ideal $I(G)$ satisfies $\mathrm v(I(G))\leq\operatorname{im}(G)
\leq\operatorname{reg}(R/I(G))$,
where $R$ is the corresponding polynomial ring and
$\operatorname{im}(G)$ is the induced matching number; see \cite{jaramillo2021v,saha2022v}.
The asymptotic behavior of v-numbers has been studied in
\cite{conca2024note,ficarra2025asymptotic,kumar2025slope}.

The v-number was introduced in the study of minimum distance
functions associated with projective Reed--Muller-type codes \cite{cooper2020generalized}.
Regardless of this interpretation, $\delta_I$ is defined for any graded ideal in terms of the Hilbert--Samuel multiplicity of $S/I$; see~\cite{nunez2017footprint}. The local notion of $\mathrm{v}$-numbers also generalizes the degree of a point in a finite set of projective points introduced in~\cite{geramita1993cayley}.

In this work, we study the $\mathrm{v}$-number of a binomial edge ideal $J_G$ at its minimal primes. This invariant has been studied for several classes of graphs, including closed graphs, paths, cycles, binary trees, and crown graphs; see \cite{dey2024v,ambhore2024v,jaramillo2024connected,liwski2026v,kumar2025binomial}. Here, we seek a general combinatorial approach to the local $\mathrm{v}$-numbers of $J_G$.

Since the $\mathrm{v}$-number is additive over the connected components of $G$ (see~\textup{\cite[Proposition~3.1]{jaramillo2024connected}} and~\textup{\cite[Corollary~3.12]{ambhore2024v}}), we may assume that $G$ is connected. In this case, the minimal primes of $J_G$ correspond to the cutsets of $G$; see Definition~\ref{minimal cuts} and Theorem~\ref{Theorem minimal primes}. For a cutset $\mathcal{S}\subseteq[n]$, let $P_{\mathcal{S}}$ denote the corresponding minimal prime, as in Definition~\ref{PS}. If

$$
G_{[n]\setminus\mathcal{S}}=G_1\amalg\cdots\amalg G_k,
$$

where $G_1,\ldots,G_k$ are the connected components, we ask:

\begin{question}
How much of the local $\mathrm{v}$-number $\VPJ$ can be determined from the combinatorial data of $\mathcal{S}$ and $G_1,\ldots,G_k$?
\end{question}

Our approach is to associate to $\mS$ two natural combinatorial invariants, $\lvs$ and $\uvs$; see Definition~\ref{def: main invariants}. 
The lower invariant is an ordinary transversal number, whereas
the upper invariant is a weighted transversal number in which
singletons have weight one and pairs have weight two.
The first main result shows that these combinatorial quantities provide lower and upper bounds for the local $\mathrm{v}$-number.

\begin{theoremA}\label{thm:A}
Let $\mS\subseteq [n]$ be a cutset. Then
\[
\lvs \leq \VPJ \leq \uvs \tag{Theorem~\ref{thm: inequalities with lvs and uvs}}.
\]
\end{theoremA}

The upper bound in Theorem~\ref{thm:A} was established in \cite{liwski2026v}, building on earlier upper bounds for the local $\mathrm{v}$-numbers of binomial edge ideals obtained in \cite{ambhore2024v,dey2024v,jaramillo2024connected}. In particular, the bound $\uvs$ recovers several known results: for $\mS=\emptyset$, it coincides with the exact value of $\VPJ$ obtained in \cite{ambhore2024v,jaramillo2024connected}; for cycles and binary trees, it agrees with the upper bounds of \cite{dey2024v}; and for Cohen--Macaulay closed graphs and paths, it recovers the precise local $\mathrm{v}$-numbers computed in \cite{dey2024v}.

The central contribution of this paper is a combinatorial framework
for lower bounds on the local $\mathrm{v}$-numbers of binomial edge
ideals. Specifically, the left-hand inequality in
Theorem~\ref{thm:A} gives a new combinatorial lower bound for $\VPJ$
by restricting attention to cutsets that do not merge distinct
connected components of $G_{[n]\setminus\mathcal S}$.
Comparing this lower bound with the upper bound from \cite{liwski2026v} leads
to a characterization of when the two bounds coincide, together
with a sufficient condition for equality. In these cases, the
local $\mathrm{v}$-number is determined entirely by the
combinatorics of the graph.

Having established Theorem~\ref{thm:A}, we turn to the combinatorial invariants $\lvs$ and $\uvs$ themselves. Their definitions in Definition~\ref{def: main invariants} are given in terms of collections of transversals, and therefore do not yet provide an explicit description in terms of the graph. Finding these transversals is itself a nontrivial combinatorial problem. Under~\eqref{hyp:H}, Theorems \ref{thm: explicit formula of lvs} and~\ref{thm: explicit formula of uvs} replace the transversal conditions
with explicit graph-theoretic descriptions involving connected domination
and the way other cutsets merge the components. 

The bounds in Section~\ref{section 3} hold for every cutset $\mathcal S$.
For the graph-theoretic descriptions in Sections~\ref{sec: explicit formulas}-\ref{sec: gap between bounds}, we assume
\begin{equation}\tag{H}\label{hyp:H}
\text{every connected component of }
G_{[n]\setminus\mathcal S}\text{ is noncomplete}. 
\end{equation}
The tree applications in Section~\ref{subsec: trees} also use this hypothesis.
The cycle computations in Section~\ref{subsec: cycles} are obtained directly from
the transversal definitions and do not require~\eqref{hyp:H}. Throughout the remainder of this introduction, we assume
hypothesis~\eqref{hyp:H}.

The descriptions in Section~\ref{sec: explicit formulas} make the two bounds more accessible and allow us to study when they coincide. In Section~\ref{sec: when the bounds coincide}, we study precisely this question. We say that a cutset $\mathcal{S}$ is \emph{v-combinatorial} if $\lvs=\uvs$. In this case, Theorem~\ref{thm:A} immediately implies

\[
\VPJ=\lvs=\uvs,
\]

so the local $\mathrm{v}$-number is determined by the combinatorial data associated to $\mathcal{S}$.

Proposition~\ref{prop: when S is v-combinatorial} gives an explicit characterization of the $\mathrm{v}$-combinatorial cutsets. Although the characterization is somewhat technical, it leads to a simple and useful sufficient condition in terms of the connected domination numbers of the components of $G_{[n]\setminus\mathcal{S}}$. Recall the connected domination number $\gamma_c(G)$ from Definition~\ref{definition connected domination}.

\begin{theoremB} Let $\mathcal S$ be a cutset, and let
$G_1,\ldots,G_k$ be the connected components of
$G_{[n]\setminus\mathcal S}$. If $\gamma_{c}(G_{i}) \geq k-1$ for all $i \in [k]$, then $\mathcal{S}$ is $\mathrm{v}$-combinatorial and

\[
\VPJ = \lvs = \uvs.
\tag{Theorem~\ref{thm: all connected domination at least k-1}}
\]

In particular, applying this result for $k=1,2,3$, we obtain:
\begin{enumerate}[label={\rm(\roman*)}]
\item The empty set $\emptyset$ is $\mathrm{v}$-combinatorial, and

\[
\textup{v}_{P_{\emptyset}}(J_G) = \gamma_{c}(G).
\tag{Theorem~\ref{thm: k=1}}
\]

\item If $k=2$, then $\mathcal{S}$ is $\mathrm{v}$-combinatorial. Let $V_{1}$ and $V_{2}$ denote the vertex sets of the connected components of the induced subgraph on $[n] \setminus \mathcal{S}$. Then

\[
\VPJ = \gamma_{c}(G,\mathcal{S})
:=
\min \left\{ |A| : A \in \mathcal{D}_{c}(G,\mathcal{S}) \right\},
\tag{Theorem~\ref{thm: k=2}}
\]

where

\[
\mathcal{D}_{c}(G,\mathcal{S})
=
\left\{
A \subseteq [n]\setminus \mathcal{S}
\;\middle|\;
A \cap V_{i} \in \mathcal{D}_{c}\big(G_{V_{i}\cup \mathcal{S}}\big)
\text{ for } i=1,2
\right\}.
\]

\item If $k=3$ and none of $G_1,G_2,G_3$ is a cone, then $\mathcal{S}$ is $\mathrm{v}$-combinatorial and

\[
\VPJ = \lvs = \uvs.
\tag{Theorem~\ref{thm:k=3}}
\]

Here, a \emph{cone} is a graph containing a vertex adjacent to all other vertices.
\end{enumerate}
\end{theoremB}

The cases $k=1$ and $k=2$ recover, under hypothesis~\eqref{hyp:H},
 known formulas for the local $\mathrm v$-numbers
\cite{ambhore2024v,jaramillo2024connected,liwski2026v}. The case $k=3$ gives
a further sufficient condition for equality of the two bounds.

We also quantify how far the two bounds can be apart.
The following estimate depends on the number of components and
their connected domination numbers.


\begin{theoremC}
Let $\mathcal{S}$ be a cutset, and let $G_{1}, \ldots, G_{k}$ denote the connected components of $G_{[n]\setminus \mathcal{S}}$. Then
\[
\uvs \leq \lvs + \sum_{i=1}^{k} \max\{\, k-1 - \gamma_{c}(G_{i}), \, 0 \,\}.
\]
Consequently,
\[
\VPJ \in [\, \lvs,\; \lvs + \sum_{i=1}^{k} \max\{\, k-1 - \gamma_{c}(G_{i}), \, 0 \,\}].
\tag{Theorem~\ref{thm: gap}}
\]
In particular, if $k\geq 2$, $\VPJ \in [\, \lvs,\; \lvs + (k-2)k]$.
\end{theoremC}

The paper is organized as follows. Section~\ref{preliminaries} recalls the required background. Section~\ref{section 3} establishes
the two combinatorial bounds, and Section~\ref{sec: explicit formulas} gives their
graph-theoretic descriptions under~\eqref{hyp:H}. Sections~\ref{sec: when the bounds coincide} and~\ref{sec: gap between bounds} study
equality and the gap between the bounds. Section~\ref{sec: examples} applies the
results to cycles and trees. Finally, Section~\ref{regularity} discusses the relation between the combinatorial upper bounds and the Castelnuovo--Mumford regularity, and formulates a conjecture and an open question concerning their minimum over all cutsets.

\medskip

\noindent\textbf{\bf Acknowledgement.}  
The author thanks Manohar Kumar for his many valuable comments and suggestions on this manuscript.
The author was supported by the PhD fellowship 1126125N and partially funded by the FWO grants G0F5921N (Odysseus), G023721N, and the KU Leuven grant iBOF/23/064.

\section{Preliminaries}\label{preliminaries}

In this section, we introduce the notation and background results on binomial edge ideals, graphs, and matroids that will be used throughout the paper. We use the notation established in the introduction and work over $\CC$. 

For a graph $H$ and a subset $U\subseteq V(H)$, we denote by
$H_U$ the induced subgraph on $U$. For $v\in V(H)$, let
$N_H(v)=\{u\in V(H):\{u,v\}\in E(H)\}$ denote its open neighborhood.
All graphs considered in this paper are finite and simple.

\subsection{Binomial edge ideals}

Let $G$ be a simple graph on $[n]$, and let $J_G$ denote its binomial edge ideal, as defined in Definition~\ref{definition binomial edge}. It was shown in \cite{herzog2010binomial} that binomial edge ideals are radical and that their associated primes admit a combinatorial description.

\begin{definition}\label{PS}
Let $\mathcal{S}\subseteq [n]$, and let $G_1,\ldots,G_{c(\mathcal{S})}$ denote the connected components of the induced subgraph $G_{[n]\setminus\mathcal{S}}$. For each $i$, let $\widetilde{G}_i$ denote the complete graph on $V(G_i)$. We define

$$
P_{\mathcal{S}}
:=
(x_i,y_i:i\in\mathcal{S})
+
(J_{\widetilde{G}_1},\ldots,J_{\widetilde{G}_{c(\mathcal{S})}}).
$$

\end{definition}

The ideal $P_{\mathcal{S}}$ is prime, and we have the following (redundant) decomposition:

\begin{theorem}\textup{\cite[Theorem~3.2]{herzog2010binomial}}\label{theorem decomposition}
Let $G$ be a simple graph on $[n]$. Then

$$
J_G=\bigcap_{\mathcal{S}\subseteq[n]}P_{\mathcal{S}}.
$$

\end{theorem}

Throughout the paper, we assume that $G$ is connected.

\begin{definition}\label{minimal cuts}
A vertex $i\in V(H)$ is called a \emph{cut point} of a graph $H$ if the number of connected components of $H_{V(H)\setminus{i}}$ is greater than that of $H$. A subset $\mathcal{S}\subsetneq[n]$ is called a \emph{cutset} of $G$ if, for every $i\in\mathcal{S}$, the vertex $i$ is a cut point of the induced subgraph $G_{([n]\setminus\mathcal{S})\cup\{i\}}$. Equivalently, $\mS$ is a cutset if every $i\in \mS$ has neighbors in at least two distinct components of $G_{[n]\setminus \mS}$.
\end{definition}

Note that, by definition, the empty set is a cutset.

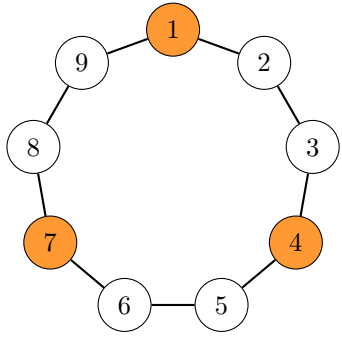
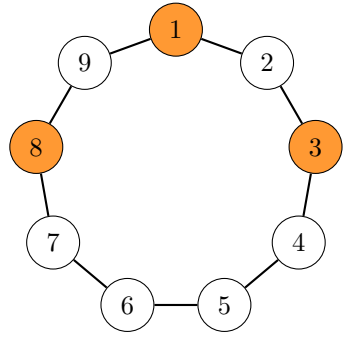
\begin{figure}[ht]
    \centering

    \begin{subfigure}[b]{0.45\textwidth}
        \centering
        \begin{tikzpicture}[
            scale=1.25,
            every node/.style={circle, draw, minimum size=7mm,
                               inner sep=0pt, font=\small},
            cut/.style={fill=orange!80},
            edge/.style={draw, thick}
        ]

            \node[cut] (v1) at (90:1.5)  {1};
            \node       (v2) at (50:1.5)  {2};
            \node       (v3) at (10:1.5)  {3};
            \node[cut] (v4) at (-30:1.5) {4};
            \node       (v5) at (-70:1.5) {5};
            \node       (v6) at (-110:1.5){6};
            \node[cut] (v7) at (-150:1.5){7};
            \node       (v8) at (170:1.5) {8};
            \node       (v9) at (130:1.5) {9};

            \draw[edge] (v1)--(v2)--(v3)--(v4)--(v5)
                        --(v6)--(v7)--(v8)--(v9)--(v1);

        \end{tikzpicture}
        \caption{$C_9$ with cutset $\{1,4,7\}$.}
    \end{subfigure}
    \hfill
    \begin{subfigure}[b]{0.45\textwidth}
        \centering
        \begin{tikzpicture}[
            scale=1.25,
            every node/.style={circle, draw, minimum size=7mm,
                               inner sep=0pt, font=\small},
            cut/.style={fill=orange!80},
            edge/.style={draw, thick}
        ]

            \node[cut] (v1) at (90:1.5)  {1};
            \node       (v2) at (50:1.5)  {2};
            \node[cut] (v3) at (10:1.5)  {3};
            \node       (v4) at (-30:1.5) {4};
            \node       (v5) at (-70:1.5) {5};
            \node       (v6) at (-110:1.5){6};
            \node       (v7) at (-150:1.5){7};
            \node[cut] (v8) at (170:1.5) {8};
            \node       (v9) at (130:1.5) {9};

            \draw[edge] (v1)--(v2)--(v3)--(v4)--(v5)
                        --(v6)--(v7)--(v8)--(v9)--(v1);

        \end{tikzpicture}
        \caption{$C_9$ with cutset $\{1,3,8\}$.}
    \end{subfigure}

    \caption{Two cutsets of $C_9$.}
    \label{fig:c9-cutsets}
\end{figure}

\begin{example}\label{ex: Cn 1}
Let $C_n$ be the cycle graph on $[n]$, with edge set 
\[E(C_n)=\{\{i,i+1\}:1\leq i<n\}\cup\{\{n,1\}\},\qquad n\geq3.\]
The cutsets of $C_n$ consist of the empty set together with all subsets
$\mS\subseteq[n]$ of cardinality at least two such that $\{i,j\}\notin E(C_n)$ for all $i,j\in \mS$.
See Figure~\ref{fig:c9-cutsets} for examples.
\end{example}

The minimal primes of $J_G$ are precisely those ideals $P_{\mathcal S}$ for which $\mathcal S$ is a cutset of $G$.

\begin{theorem}\textup{\cite[Corollary~3.9]{herzog2010binomial}}\label{Theorem minimal primes}
Let $G$ be a connected graph on $[n]$, and let $\mathcal{S}\subseteq[n]$. Then $P_{\mathcal{S}}$ is a minimal prime of $J_G$ if and only if $\mathcal{S}$ is a cutset of $G$.
\end{theorem}

\subsection{Matroids and connected domination}\label{matroids}

We use matroids of rank at most two. We recall their definitions
and refer to \cite{Oxley} for general matroid theory.


\begin{definition}
A matroid $M$ of rank at most two on $[n]$ is specified by a family $\mathcal{I}(M)\subseteq \{I\subseteq [n]:\size{I}\leq 2\}$ of  independent sets. This family contains $\emptyset$, is closed under taking subsets, and satisfies the augmentation condition: if $I_1,I_2\in \mathcal{I}(M)$ with $\size{I_1}<|I_2|$, then \(    I_1\cup\{e\}\in\mathcal{I}(M)
   \) for some $e\in I_2\setminus I_1$.

\end{definition}

\begin{definition}\label{def:dependant}
For a matroid $M$ on $[n]$, write $\mathcal{D}(M)=\{B\subseteq[n]: B\notin \mathcal{I}(M)\}$ for its dependent sets. An element $x\in[n]$ is a loop when $\{x\} \in\mathcal{D}(M)$. Distinct nonloops $x,y$ are parallel (or form a double point) when $\{x,y\}\in\mathcal{D}(M)$.
\end{definition}


Let $G$ be a simple graph on $[n]$. To each subset $\mathcal{S}\subseteq[n]$, we associate a matroid $M(\mathcal{S})$ of rank at most two, defined as follows.

\begin{definition}\label{definition MS}
Let $\mathcal{S}\subseteq[n]$, and let $G_1,\ldots,G_{c(\mathcal{S})}$ denote the connected components of the induced subgraph $G_{[n]\setminus\mathcal{S}}$. We define a matroid $M(\mathcal{S})$ on $[n]$ via its collection of dependent sets:
\[
\mathcal D(M(\mathcal S))
=
\{B\subseteq[n]:B\cap\mathcal S\neq\emptyset\}
\cup
\bigcup_{i\in[c(\mathcal S)]}
\{B\subseteq[n]:|B\cap V(G_i)|\geq2\}
\cup
\{B\subseteq[n]:|B|\geq3\}.
\]
Thus, the vertices of $\mathcal S$ are loops, and the nonempty
sets $V(G_1),\ldots,V(G_{c(\mathcal S)})$ are the parallel classes
of the remaining elements. In particular, an independent pair
consists of vertices from two distinct components.
\end{definition}

\begin{example}\label{ex: Cn 2}
Let $\mS\subseteq[n]$ be a nonempty cutset of $C_n$, and let
$I_1,\ldots,I_k$ be the disjoint intervals whose union is
$[n]\setminus\mS$. Then
\[
\mathcal{D}(M(\mathcal{S}))
=
\{B\subseteq[n]:B\cap\mathcal{S}\neq\emptyset\}
\cup
\bigcup_{j=1}^{k}
\{B\subseteq[n]:|B\cap I_j|\geq2\}
\cup
\{B\subseteq[n]:|B|\geq3\}.
\]
\end{example}

These matroids will play a central role in our study of the local $\mathrm{v}$-numbers of binomial edge ideals. 
We will also use the following graph-theoretic invariant.

\begin{definition}\label{definition connected domination}
For a connected graph $G$ on $[n]$, let
\[\mathcal{D}_c(G):=\{ B\subseteq [n]: G_{B} \text{ is connected and } N_{G}(v)\cap B\neq \emptyset \text{ for every } v\in [n]\setminus B\}.\]
Its members are called the 
\emph{connected dominating sets} of $G$. The \emph{connected domination number} of $G$ is
$$
\gamma_c(G):=\min\{|B|:B\in\mathcal{D}_c(G)\}.
$$

\end{definition}

\section{\texorpdfstring{Combinatorial bounds $\lvs$ and $\uvs$}
{Combinatorial bounds lvs and uvs}}
\label{section 3}

In this section, we introduce the combinatorial invariants
$\lvs$ and $\uvs$ associated to a cutset $\mathcal{S}$ and show that
they provide, respectively, a lower and an upper bound for the local
$\mathrm{v}$-number $\VPJ$. Both invariants are defined in terms of transversals of suitable hypergraphs.
We begin by recalling the definition of a transversal and introducing
the necessary notation, following \cite{liwski2026v}.

\begin{definition}\label{transveral definition}
Let $\Delta\subseteq 2^E$ be a hypergraph on a finite ground set $E$.
A subset $A\subseteq E$ is a transversal of $\Delta$ if
$A\cap e\neq\emptyset$ for every $e\in\Delta$.
If $\Delta=\emptyset$, every subset of $E$, including the empty
set, is a transversal.
\end{definition}

\begin{definition}\label{important notation}
We use the following notation. 

\begin{enumerate}
    \item[{\rm (i)}] For $k\geq1$, let
\[
{\min}_k(G)=\{\mathcal S:\mathcal S\text{ is a cutset of }G
\text{ and }G_{[n]\setminus\mathcal S}
\text{ has }k\text{ connected components}\}.
\]
Since $G$ is connected, ${\min}_1(G)=\{\emptyset\}$.
We put $\min(G)=\cup_{k\geq1}{\min}_k(G)$.

    \item[{\rm (ii)}]
    For a fixed $\mathcal S\in\min(G)$, put
\[
E_{\mathcal S}=
\left(\binom{[n]}1\cup\binom{[n]}2\right)
\setminus\mathcal D(M(\mathcal S)).
\]
Define the hypergraph $\Delta_{\mS}$ on $E_{\mathcal S}$ by
\[
\Delta_{\mathcal S}=
\{\mathcal D(M(\mathcal S'))\setminus\mathcal D(M(\mathcal S)):
\mathcal S'\in\min(G)\setminus\{\mathcal S\}\}.
\]
We denote by $\mathcal T(\mathcal S)$ its family of transversals.


    \item[{\rm (iii)}]\label{item 3} Each $A \in \mathcal{T}(\mathcal{S})$ can be written as a disjoint union $A = A_1 \amalg A_2$, with \[A_1 \subseteq \textstyle \binom{[n]}{1}\quad  \text{and} \quad A_2 \subseteq \binom{[n]}{2},\]
    and we put
\[
\operatorname{supp}(A)=\bigcup_{a\in A}a.
\]
We identify $A_1$ with the subset $\{i\in [n]: \{i\}\in A_1\}$ of $[n]$.

    \item[{\rm (iv)}] Given $C, D \subseteq [n]$, put \[g_{C,D} := \prod_{k \in C} x_k \prod_{k \in D} y_k.\]

    \item[{\rm (v)}]\label{item v} Given $A \in \mathcal{T}(\mathcal{S})$ and a partition $A_1 = C \amalg D$, define \[g_{A,C,D} := g_{C,D} \cdot \prod_{\{i,j\} \in A_2} f_{i,j},\]
    where $f_{i,j}=x_{i}y_{j}-x_{j}y_{i}$.
\end{enumerate}
\end{definition}

Throughout this section, we fix a graph $G$ on $[n]$ and a cutset
$\mathcal{S}\in{\min}_k(G)$. Let
$G_1,\ldots,G_k$ be the connected components of the induced subgraph
on $[n]\setminus\mathcal{S}$, and denote their respective vertex sets
by $V_1,\ldots,V_k$.


\begin{definition}\normalfont \label{Tprima}
With respect to the fixed cutset $\mathcal{S}\in {\min}_{k}(G)$, we introduce the following terminology:
\begin{itemize}
\item A subset $\mathcal S'\subseteq[n]$ is nice with respect
to $\mathcal S$ if the vertex set of every connected
component of $G_{[n]\setminus\mathcal S'}$ intersects
at most one of $V_1,\ldots,V_k$. Observe that if $\mathcal S'\neq\mathcal S$ is nice with respect to
$\mathcal S$, then $\mathcal S'\not\subseteq\mathcal S$.

\item Let $\Delta_{S}'$ denote the hypergraph on $[n]\setminus \mS$ with collection of edges
\begin{equation}\label{delta}
    \left\{ \mathcal{S}' \setminus \mathcal{S} : \mathcal{S}' \in \min(G)\setminus \{\mathcal{S}\}, \ \text{and $\mathcal{S}'$ is nice with respect to $\mathcal{S}$} \right\},
\end{equation}
and let $\mathcal{T}'(\mathcal{S})$ be the set of all transversals of $\Delta_{\mathcal{S}}'$. 
\end{itemize}
\end{definition}

\begin{example}
The cutsets $\{1,4,7\}$ and $\{1,3,8\}$ of $C_9$, shown in
Figure~\ref{fig:c9-cutsets}, are nice with respect to each other.
\end{example}

We now introduce the combinatorial invariants $\lvs$ and $\uvs$, defined in terms of the collections of transversals $\mathcal{T}(\mathcal{S})$ and  $\mathcal{T}'(\mathcal{S})$.

\begin{definition}\label{def: main invariants}
The invariants $\lvs$ and $\uvs$ are defined by
\[
\lvs := \min \{\, |A| \mid A \in \mathcal{T}'(\mathcal{S}) \,\},
\]
and
\[
\uvs := \min \{\, |A_{1}| + 2|A_{2}| \mid A \in \mathcal{T}(\mathcal{S}) \,\}.
\]
\end{definition}

The main result of this section is that these two quantities bound the
local $\mathrm{v}$-number from opposite sides.

\begin{theorem}\label{thm: inequalities with lvs and uvs}
For every cutset $\mS$ of $G$, the local \textup{v}-number satisfies
\begin{equation}\label{inequality with lvs and uvs}
\lvs \leq \VPJ \leq \uvs .
\end{equation}
\end{theorem}

To prove the theorem, we recall the ideal
$J_{\mathcal T(\mathcal S)}$ introduced in
\cite{liwski2026v} and define a second ideal
$J_{\mathcal T'(\mathcal S)}$ associated to
the transversals $\mathcal T'(\mathcal S)$.

\begin{definition}\label{notation JS}
We denote
\[
J_{\mathcal{T}(\mathcal{S})}
:=
\left(
g_{A,C,D}
\;\middle|\;
A \in \mathcal{T}(\mathcal{S}), \; C \amalg D = A_1
\right)
\]
and
\begin{equation}\label{equ: pJS}
J_{\mathcal{T}'(\mathcal{S})}
:=
\left(
g_{C,D}
\;\middle|\;
C \amalg D \in \mathcal{T}'(\mathcal{S})
\right).
\end{equation}
Note that $\pJS$ is a monomial ideal.
\end{definition}


The following lemma relates the two ideals and will be used in the
proof of Theorem~\ref{thm: inequalities with lvs and uvs}.

\begin{lemma}\label{lem: inc lem}
We have $\JS \subseteq \pJS$.
\end{lemma}

\begin{proof}

Let $A\in\mathcal T(\mathcal S)$, and let
$\mathcal S'\in\min(G)\setminus\{\mathcal S\}$ be nice with
respect to $\mathcal S$. Every element of
$\mathcal D(M(\mathcal S'))\setminus\mathcal D(M(\mathcal S))$
contains a vertex of $\mathcal S'\setminus\mathcal S$.
Indeed, a singleton in this difference is such a vertex; and a
pair avoiding $\mathcal S'$ could be dependent in $M(\mathcal S')$ and independent in $M(\mathcal S)$
only by joining two original components, contrary to niceness.
Since $A$ intersects $\mathcal D(M(\mathcal S'))\setminus\mathcal D(M(\mathcal S))$, its support intersects
$\mathcal S'\setminus\mathcal S$. Thus
$\operatorname{supp}(A)\in\mathcal T'(\mathcal S)$.

Every monomial occurring in $g_{A,C,D}$ contains a variable
indexed by each vertex of $\operatorname{supp}(A)$.
Choosing one such variable for each vertex gives a generator of
$J_{\mathcal T'(\mathcal S)}$ dividing that monomial.
Hence $g_{A,C,D}\in J_{\mathcal T'(\mathcal S)}$.
This holds for every generator of $J_{\mathcal T(\mathcal S)}$,
and therefore
$J_{\mathcal T(\mathcal S)}\subseteq J_{\mathcal T'(\mathcal S)}$.

\end{proof}

To obtain the lower bound, we use the
relationship between the transversal ideals and the colon ideal of
$P_{\mathcal{S}}$, given by the following theorem from~\cite{liwski2026v}.

\begin{theorem}\textup{\cite[Theorem~3.3]{liwski2026v}}\label{theorem transversal}
With the notation above, we have:
\begin{equation}\label{equality of ideals}(J_{G}:P_{\mathcal{S}})=\sqrt{J_{G}+\JS}.\end{equation}
\end{theorem}


We now show that the invariants $\lvs$ and $\uvs$ provide bounds for $\VPJ$.

\begin{proof}[Proof of Theorem~\ref{thm: inequalities with lvs and uvs}]
The upper bound is \textup{\cite[Theorem~3.5]{liwski2026v}}. We prove the lower bound.
Use the lexicographic order
$x_1>\cdots>x_n>y_1>\cdots>y_n$.
By \textup{\cite[Theorem~2.1]{herzog2010binomial}}, $J_G$ has the reduced Gr\"obner basis
\[
\mathcal G=\{u_\pi f_{i,j}:i<j,\ \pi\text{ is an admissible
path from }i\text{ to }j\}.
\]
Here $u_\pi$ is the squarefree monomial associated with the path
$\pi$ as in \cite{herzog2010binomial}, and the leading monomial is $u_\pi x_i y_j$.
Let $\mathcal G'$ be the squarefree monomial generating set of
$J_{\mathcal T'(\mathcal S)}$ from~\eqref{equ: pJS}.

We claim that $\mathcal G\cup\mathcal G'$ is a Gr\"obner basis
of $J_G+J_{\mathcal T'(\mathcal S)}$.
The pairs within $\mathcal G$ reduce to zero, and the pairs
within $\mathcal G'$ are zero. For a mixed pair, write
$g=u_\pi f_{i,j}$ and $h=g_{C,D}$, with
$C\amalg D\in\mathcal T'(\mathcal S)$, and put
$L=\operatorname{lcm}(u_\pi x_i y_j,h)$.
Then
\[
S(g,h)=-\frac{L}{x_i y_j}x_jy_i.
\]
The monomial on the right contains a variable indexed by every
vertex of $C\cup D$: removing $x_i$ is compensated by $y_i$,
and removing $y_j$ by $x_j$.
Choose one such variable for each vertex of $C\cup D$.
Their product is a generator of $J_{\mathcal T'(\mathcal S)}$
and divides the monomial. Thus the mixed pair reduces to zero. Consequently,
\[
\operatorname{in}_{<}(J_G+J_{\mathcal T'(\mathcal S)})
=\operatorname{in}_{<}(J_G)+J_{\mathcal T'(\mathcal S)}.
\]
All leading monomials in this Gr\"obner basis are squarefree,
so $J_G+J_{\mathcal T'(\mathcal S)}$ is radical.

Let $f$ be homogeneous with $(J_G:f)=P_{\mathcal S}$.
Replace $f$ by its remainder modulo the homogeneous basis
$\mathcal G$. The remainder is nonzero, homogeneous of the same
degree, and has the same colon ideal, because it differs from
$f$ by an element of $J_G$. We may therefore assume
$\operatorname{in}_{<}(f)\notin\operatorname{in}_{<}(J_G)$.
By Theorem~\ref{theorem transversal}, Lemma~\ref{lem: inc lem}, and the radicality just established,
\[
f\in(J_G:P_{\mathcal S})
=\sqrt{J_G+J_{\mathcal T(\mathcal S)}}
\subseteq
\sqrt{J_G+J_{\mathcal T'(\mathcal S)}}
=J_G+J_{\mathcal T'(\mathcal S)}.
\]
It follows that
$\operatorname{in}_{<}(f)\in J_{\mathcal T'(\mathcal S)}$.
Every generator of this monomial ideal has degree at least
$\mathrm v^-_{\mathcal S}(G)$, and hence
\[
\deg(f)=\deg(\operatorname{in}_{<}(f))
\geq\mathrm v^-_{\mathcal S}(G).
\]
Taking the minimum over all such $f$ proves the lower bound.

\end{proof}

\begin{remark}
By \textup{\cite[Theorem~3.5]{liwski2026v}}, the equality $\VPJ=\uvs$ holds if $J_{G}+\JS$ is a radical ideal.
\end{remark}

\section{\texorpdfstring{Explicit formulas for the bounds $\lvs$ and $\uvs$}
{Explicit formulas for the bounds lvs and uvs}}\label{sec: explicit formulas}

Fix a cutset $\mathcal S$ and write
$G_{[n]\setminus\mathcal S}=G_1\sqcup\cdots\sqcup G_k$,
with $V_i=V(G_i)$. Throughout this section, except where
explicitly stated otherwise, assume~\eqref{hyp:H}.
We describe the lower invariant using connected domination
and then describe the additional conditions for the upper bound.


We begin with the lower bound. Recall the notion of $\mathcal{D}_{c}(G)$ from Definition~\ref{definition connected domination}.

\begin{definition}\label{def: DC-}
Let $\DCS$ denote the collection of all subsets $A\subseteq [n]\setminus \mS$ satisfying the following conditions:
\begin{enumerate}[label={\rm(\roman*)}]
\item\label{item:1} $A\cap V_{i}\in \mathcal{D}_{c}(G_{i})$ for all $i\in [k]$.
\item\label{item:2} For every $s\in \mS$, there exist distinct $i,j\in [k]$ such that
\[
N_{G}(s)\cap A\cap V_{i}\neq \emptyset 
\quad \text{and} \quad 
N_{G}(s)\cap A\cap V_{j}\neq \emptyset.
\]
\end{enumerate}
\end{definition}

The next theorem identifies $\mathcal{T}'(\mS)$ with $\DCS$.

\begin{theorem}\label{thm: explicit formula of lvs}
Let $\mathcal S$ be a cutset satisfying~\eqref{hyp:H}.
For $A\subseteq[n]\setminus\mathcal S$,
$A\in\mathcal T'(\mathcal S)$ if and only if
$A\in\mathcal D_c^-(\mathcal S)$.
In particular,
\[
\mathrm v^-_{\mathcal S}(G)
=\min\{|A|:A\in\mathcal D_c^-(\mathcal S)\}.
\]
\end{theorem}

\begin{proof}
We first prove that if $A$ satisfies conditions~\ref{item:1} and~\ref{item:2}, then $A\in \pTS$. It suffices to show that
\begin{equation}\label{nonempty intersection}
A\cap \mS'\neq \emptyset
\end{equation}
for every $\mS'\in \min(G)\setminus \{\mS\}$ that is nice with respect to $\mS$. Assume, for contradiction, that $A\subseteq [n]\setminus \mS'$. We first consider the case $\mS\subsetneq \mS'$. In this case,
\[
[n]\setminus \mS'=(V_{1}\setminus \mS')\amalg \cdots \amalg (V_{k}\setminus \mS').
\]
Since $A\subseteq [n]\setminus \mS'$, we have $A\cap V_{i}\subseteq V_{i}\setminus \mS'$ for all $i$. By condition~\ref{item:1}, the set $A\cap V_{i}$ is connected and dominating in $G_i$, and therefore the graph $G_{V_{i}\setminus \mS'}$ is connected for each $i\in [k]$. Consequently, the graph $G_{[n]\setminus \mS'}$ has exactly $k$ connected components, namely
\[
G_{V_{1}\setminus \mS'},\ldots,G_{V_{k}\setminus \mS'}.
\]

Since $\mS'\in \min(G)$, every vertex $s\in \mS'$ must be a cut point of the graph $G_{([n]\setminus \mS')\cup \{s\}}$. Let $s\in \mS'\setminus \mS$ and assume without loss of generality that $s\in V_{i}$ for some $i\in [k]$. As $s$ is adjacent only to vertices of the connected component $G_{V_{i}\setminus \mS'}$ of $G_{[n]\setminus \mS'}$, it cannot be a cut point, which is a contradiction. Hence~\eqref{nonempty intersection} holds whenever $\mS\subsetneq \mS'$.

Next suppose that $\mS\not\subset \mS'$, and again assume for contradiction that $A\subseteq [n]\setminus \mS'$. Let $s\in \mS\setminus \mS'$. By condition~\ref{item:2}, there exist distinct $i,j\in [k]$ such that
\[
N_{G}(s)\cap A\cap V_{i}\neq \emptyset 
\quad \text{and} \quad 
N_{G}(s)\cap A\cap V_{j}\neq \emptyset.
\]
Choose $i_{1}\in A\cap V_{i}$ and $i_{2}\in A\cap V_{j}$ such that $\{i_{1},s\},\{i_{2},s\}\in E(G)$. Since $A\subseteq [n]\setminus \mS'$, we have $i_{1}\in V_{i}\setminus \mS'$ and $i_{2}\in V_{j}\setminus \mS'$. Because both vertices are adjacent to $s$, they lie in the same connected component of $G_{[n]\setminus \mS'}$, contradicting the assumption that $\mS'$ is nice with respect to $\mS$. Therefore~\eqref{nonempty intersection} holds, and hence $A\in \pTS$.

Conversely, suppose that $A\in \pTS$. We show that $A$ satisfies conditions~\ref{item:1} and~\ref{item:2}. 

We first verify~\ref{item:1}. We fix $i\in [k]$ and will show that $A\cap V_{i}\in \mathcal{D}_{c}(G_i)$. By \textup{\cite[Lemma~3.7]{liwski2026v}}, this is equivalent to proving that $A\cap V_{i}$ is a transversal of the family
\[
\Delta:=\{T\subset V_{i}: T\in \min(G_{i})\setminus \{\emptyset\}\}.
\]
The hypothesis~\eqref{hyp:H} is used here to ensure that the lemma applies.
Let $T\in \Delta$ with $T\in {\min}_{r}(G_i)$. We must show that $A\cap T\neq \emptyset$. Let $W_{1},\ldots,W_{r}$ denote the vertex sets of the $r$ connected components of $G_{V_{i}\setminus T}$. Define $B:=\mathcal{S}\cup T$. Then $G_{[n]\setminus B}$ has exactly $r+k-1$ connected components with vertex sets
\begin{equation}\label{vertex sets}
V_1,\ldots,V_{i-1},V_{i+1},\ldots,V_k,W_1,\ldots,W_r.
\end{equation}

By \textup{\cite[Lemma~4.2]{liwski2026v}}, there exists a subset $B'\subseteq B$ such that $B'\in \min(G)$ and each of the sets in~\eqref{vertex sets} is contained in a distinct connected component of $G_{[n]\setminus B'}$.

We first show that $T\subseteq B'$. Suppose, for contradiction, that there exists $j\in T\setminus B'$. Since $T\in \min(G_i)$, the vertex $j$ is a cut point in $G_{(V_{i}\setminus T)\cup \{j\}}$, and therefore it is adjacent to vertices in at least two of the components
\[
W_{1},\ldots,W_{r}.
\]
Hence these components would lie in the same connected component of $G_{[n]\setminus B'}$, contradicting the definition of $B'$. Thus $T\subseteq B'$.

For $j\neq i$, all vertices of $V_j$ survive in
$G_{[n]\setminus B'}$, while the surviving vertices of $V_i$
are $W_1\cup\cdots\cup W_r$.
By the choice of $B'$, these sets lie in distinct connected
components. Hence $B'$ is nice with respect to $\mathcal S$.
Moreover, $B'\neq\mathcal S$ because $T$ is nonempty and
$T\subseteq B'$.


Since $B'$ is nice and $A\in \pTS$, we obtain $A\cap B'\neq \emptyset$. As $B'\subseteq \mS\cup T$ and $A\subseteq [n]\setminus \mS$, it follows that $A\cap T\neq \emptyset$. Therefore $A\cap V_{i}\in \mathcal{D}_{c}(G_i)$.

We now verify condition~\ref{item:2}. Suppose, for contradiction, that there exists $s\in \mS$ such that $N_G(s)$ intersects at most one of the sets $V_{1}\cap A,\ldots,V_{k}\cap A$.
Put $C=[n]\setminus(A\cup\{s\})$ and $A_i=A\cap V_i$.
Each $A_i$ is nonempty and induces a connected graph.
If $s$ has neighbors in exactly one of the sets, say in $V_\ell$, the connected
components of $G_{[n]\setminus C}$ have vertex sets
$A_\ell\cup\{s\}$ and $A_i$ for $i\neq\ell$.
If $s$ has no neighbors in any of the sets, then their vertex sets are
$A_1,\ldots,A_k$ and $\{s\}$.

By \textup{\cite[Lemma~4.2]{liwski2026v}}, we can choose a cutset $C'\subseteq C$ preserving
the separation of these components. In particular, the sets
$A_1,\ldots,A_k$ remain in distinct components of
$G_{[n]\setminus C'}$.
Every vertex of $V_i\setminus C'$ is connected to $A_i$:
it either belongs to $A_i$ or has a neighbor in $A_i$,
because $A_i$ dominates $G_i$.
Therefore, no component of $G_{[n]\setminus C'}$ can intersect two distinct sets $V_i$, as this would place two of the sets $A_1,\ldots,A_k$ in the same component of $G_{[n]\setminus C'}$. Thus, $C'$
is nice with respect to $\mathcal S$.
Since $s\notin C'$ but $s\in\mathcal S$, we have
$C'\neq\mathcal S$.
The transversal condition gives $A\cap C'\neq\emptyset$,
contrary to $C'\subseteq C\subseteq[n]\setminus A$.
Thus condition~\ref{item:2} holds.
\end{proof}

Theorem~\ref{thm: explicit formula of lvs} gives a direct graph-theoretic interpretation
of $\lvs$. Using this result, we now obtain a lower bound for $\VPJ$ in terms of connected domination numbers. Recall the notion of the connected domination number from Definition~\ref{definition connected domination}.

\begin{corollary}\label{cor: sum of gamma c}
Let $\mS$ be a cutset satisfying~\eqref{hyp:H}, with components $G_1,\ldots,G_k$. We have
\[\sum_{i=1}^{k}\gamma_{c}(G_{i})\leq \lvs \leq \VPJ.\]
\end{corollary}

\begin{proof}
For every $A\in\mathcal D_c^-(\mathcal S)$, the set
$A\cap V_i$ is a connected dominating set of $G_i$.
Consequently,
\[
|A|=\sum_{i=1}^k |A\cap V_i|
\geq \sum_{i=1}^k \gamma_c(G_i).
\]
Taking the minimum over
$A\in\mathcal D_c^-(\mathcal S)$ and applying
Theorem~\ref{thm: inequalities with lvs and uvs} give the desired bounds.
\end{proof}

We now turn to the upper invariant $\uvs$. Unlike the lower invariant, whose description involves
connected domination in the components $G_i$ and
neighborhood conditions on the vertices of $\mathcal S$, the upper invariant $\uvs$
must also account for cutsets that merge different components of
$G_{[n]\setminus\mS}$. The following construction records precisely
which components are merged.

\begin{definition}
Let $\mathcal S'\in\min(G)$ be such that
$G_{V_i\setminus\mathcal S'}$ is nonempty and connected for
every $i\in[k]$. Define a partition $\mathcal P_{\mathcal S'}$
of $[k]$ by placing $r$ and $s$ in the same block if
$V_r\setminus\mathcal S'$ and $V_s\setminus\mathcal S'$
lie in the same connected component of
$G_{[n]\setminus\mathcal S'}$.
For such $\mathcal S'$, niceness with respect to $\mathcal S$
is equivalent to
$\mathcal P_{\mathcal S'}=\{\{1\},\ldots,\{k\}\}$.

Moreover, for a subset $B \subseteq \textstyle \binom{[n]}{2} \setminus \mathcal{D}(M(\mS))$, we define the simple graph $G[B]$ on vertex set $[k]$ by
\[
\{i,j\} \in E(G[B]) \quad \Longleftrightarrow \quad 
\text{there exist } v \in V_i \text{ and } w \in V_j \text{ such that } \{v,w\} \in B.
\]
\end{definition}


Thus, $\mathcal{P}_{\mS'}$ records the way in which the cutset
$\mS'$ identifies the components of $G_{[n]\setminus\mS}$, while
$G[B]$ records which pairs of those components are represented by the
two-element subsets selected in $B$.

\begin{definition}
Let $G'$ be a graph on $[k]$ and $\mathcal{P}$ a partition of $[k]$. We write
\[
G' \pitchfork \mathcal{P}
\]
to indicate that there exists a block $B \in \mathcal{P}$ such that $G'$ contains an edge $\{i,j\}$ with $i,j \in B$. In other words, $G'\pitchfork\mathcal{P}$ means that $G'$ contains
an edge entirely inside one block of $\mathcal{P}$.
\end{definition}

\begin{definition}\label{def: ECS}
Let $\ECS$ denote the collection of all subsets 
\[
A\subseteq (\tbinom{[n]}{1}\cup \tbinom{[n]}{2})\setminus \mathcal{D}(M(\mS))
\]
satisfying the following conditions:
\begin{enumerate}[label={\rm(\roman*)}]
\item\label{condi:1} $\supp(A)\in \DCS$.
\item\label{condi:2} $G[A_{2}]\pitchfork \pms$ for every $\mS'\in \min(G)$ that is not nice with respect to $\mS$ and such that $\mS'\cap\supp(A)=\emptyset$.
\end{enumerate}
\end{definition}

The partition $\pms$ in condition~\ref{condi:2} is well defined.
Indeed, if $\mathcal S'$ avoids $\operatorname{supp}(A)$,
then each $G_{V_i\setminus\mathcal S'}$ contains the
connected dominating set
$\operatorname{supp}(A)\cap V_i$ and is therefore
nonempty and connected.


The first condition ensures that the support of $A$ already satisfies
the requirements in the formula for
$\lvs$. The second condition adds precisely the extra information
needed for belonging to $\TS$: whenever a non-nice cutset $\mS'$
avoids the support of $A$, the pairs in $A_2$ must detect one pair of
components that it merges.

\begin{theorem}\label{thm: explicit formula of uvs}
Let $\mathcal S$ be a cutset satisfying~\eqref{hyp:H}. A subset 
\[
A\subseteq (\tbinom{[n]}{1}\cup \tbinom{[n]}{2})\setminus \mathcal{D}(M(\mS))
\]
belongs to $\TS$ if and only if $A\in \ECS$. In particular,
\[
\uvs=\min \{\, |A_{1}|+2|A_{2}| : A\in \ECS \,\}.
\]
\end{theorem}

\begin{proof}
Suppose first that $A\in\mathcal D_c^+(\mathcal S)$.
Let $\mathcal S'\in\min(G)\setminus\{\mathcal S\}$.
If $\mathcal S'$ is nice with respect to $\mS$,~\ref{condi:1} and Theorem~\ref{thm: explicit formula of lvs} give $\supp(A)\cap \mS'\neq \emptyset$. Hence there exists $a\in A$ such that $a\cap \mS'\neq \emptyset$. Since the elements of $\mS'$ are loops in $M(\mS')$, it follows that $a$ is dependent in $M(\mS')$. Therefore
\[
A\cap \mathcal{D}(M(\mS'))\neq \emptyset.
\]
The same conclusion holds when $\mathcal S'$ is not nice but $\supp(A)\cap \mS'\neq \emptyset$.

It remains to consider a non-nice $\mathcal S'$ avoiding
$\operatorname{supp}(A)$.
By~\ref{condi:2}, we have $G[A_{2}]\pitchfork \pms$. Hence there exists a block $B\in \pms$ and an edge $\{r,s\}\in E(G[A_{2}])$ with $r,s\in B$. This means that there exist vertices $i\in V_r\setminus \mS'$ and $j\in V_s\setminus \mS'$ lying in the same connected component of $G_{[n]\setminus \mS'}$, and vertices $i'\in V_r$ and $j'\in V_s$ such that $\{i',j'\}\in A$.

Moreover, by~\ref{condi:1}, we have that $\supp(A)\cap V_{r}$ is connected and dominating in $G_{V_{r}}$, and since $\supp(A)\subset [n]\setminus \mS'$, it follows that $G_{V_{r}\setminus \mS'}$ is connected and hence $i$ and $i'$ lie in the same connected component of $G_{[n]\setminus \mS'}$. The analogous argument shows that $j$ and $j'$ lie in the same connected component of $G_{[n]\setminus \mS'}$. Since $i$ and $j$ are also connected in this graph, it follows that $i'$ and $j'$ lie in the same connected component of $G_{[n]\setminus \mS'}$. Consequently,
\[
\{i',j'\}\in \mathcal{D}(M(\mS')),
\]
and therefore $A\cap \mathcal{D}(M(\mS'))\neq \emptyset$.

We now prove the converse. Suppose that $A\in \TS$. We first show that $A$ satisfies condition~\ref{condi:1}. Let $\mathcal S'\in\min(G)\setminus\{\mathcal S\}$
be nice with respect to $\mathcal S$. Since $A\in \TS$, we have
\[
A\cap \mathcal{D}(M(\mS'))\neq \emptyset.
\]
By niceness, every element of
$\mathcal D(M(\mathcal S'))\setminus
\mathcal D(M(\mathcal S))$
contains a vertex of $\mathcal S'\setminus\mathcal S$.
Hence it follows that $\supp(A)\cap \mS'\neq \emptyset$.
Consequently, $\supp(A)\in \mathcal{T}'(\mS)$, and Theorem~\ref{thm: explicit formula of lvs} then gives $\supp(A)\in \DCS$, establishing~\ref{condi:1}.

We now prove~\ref{condi:2}. Let $\mS'\in \min(G)$ be not nice with respect to $\mS$ and suppose $\mS'\cap\supp(A)=\emptyset$. Since $A\in \TS$, there exists
\[
a\in A\cap \mathcal{D}(M(\mS')).
\]
Because $\supp(A)\cap \mS'=\emptyset$, we have $a\cap \mS'=\emptyset$. In particular, $|a|=2$, since otherwise $a$ would be a loop of $M(\mS')$, whose loops are precisely the elements of $\mS'$, contradicting $a\cap \mS'=\emptyset$.

Write $a=\{v,w\}$ with $v\in V_i\setminus \mS'$ and $w\in V_j\setminus \mS'$. By the definition of $G[A_2]$, this implies
\[
\{i,j\}\in E(G[A_2]).
\]
Moreover, since $\{v,w\}\in \mathcal{D}(M(\mS'))$, the vertices $v$ and $w$ lie in the same connected component of $G_{[n]\setminus \mS'}$. Hence $i$ and $j$ belong to the same block of $\pms$, which implies $G[A_2]\pitchfork \pms$. This proves~\ref{condi:2} and completes the proof.
\end{proof}


The condition~\eqref{hyp:H} is used to identify transversals with
connected dominating sets in Theorem~\ref{thm: explicit formula of lvs}. It also enters
Theorem~\ref{thm: explicit formula of uvs} through that identification. This leaves the
following question.

\begin{question}
Can one still obtain explicit graph-theoretic descriptions of $\lvs$ and $\uvs$ 
for cutsets that do not satisfy condition~\eqref{hyp:H}?
\end{question}

\begin{example}
Let $\mS_1:=\{1,4,7\}$ and $\mS_2:=\{1,3,8\}$ be the cutsets of
$C_9$ shown in Figure~\ref{fig:c9-cutsets}. Neither cutset satisfies
hypothesis~\eqref{hyp:H}, so Theorems~\ref{thm: explicit formula of lvs}
and~\ref{thm: explicit formula of uvs} do not apply directly. Nevertheless,
the invariants $\textup{v}_{\mathcal{S}}^{-}(C_9)$ and
$\textup{v}_{\mathcal{S}}^{+}(C_9)$ can still be computed for these cutsets.
In fact,
\[
\textup{v}_{\mathcal{S}_1}^{-}(C_9)
=
\textup{v}_{\mathcal{S}_1}^{+}(C_9)
=
6,
\]
while
\[
\textup{v}_{\mathcal{S}_2}^{-}(C_9)=6,
\qquad
\textup{v}_{\mathcal{S}_2}^{+}(C_9)=8.
\]
Hence, by Theorem~\ref{thm: inequalities with lvs and uvs},
\[
\textup{v}_{P_{\mathcal{S}_1}}(J_{C_9})=6
\qquad\text{and}\qquad
6\leq
\textup{v}_{P_{\mathcal{S}_2}}(J_{C_9})
\leq 8.
\]
In
\S~\ref{subsec: cycles}, we explain these computations and determine
$\textup{v}_{\mathcal{S}}^{-}(C_n)$ and
$\textup{v}_{\mathcal{S}}^{+}(C_n)$ for cutsets $\mathcal{S}$
of cycle graphs with $\size{\mS}\geq 3$.
\end{example}

\section{\texorpdfstring{When the bounds $\lvs$ and $\uvs$ coincide}
{When the bounds lvs and uvs coincide}}\label{sec: when the bounds coincide}


In this section, we study the cutsets $\mathcal{S}\in\min(G)$ for which
$\lvs$ and $\uvs$ coincide. By Theorem~\ref{thm: inequalities with lvs and uvs}, such a coincidence determines
the local $\mathrm{v}$-number exactly. We first give a combinatorial
characterization of when this equality holds and then derive several sufficient
conditions under which it occurs. Throughout this section, we retain hypothesis~\eqref{hyp:H}. 

\begin{definition}\label{def: v-combinatorial}
We say that $\mS\in \min(G)$ is v-{\em combinatorial} if $\lvs=\uvs$. In that case, by Theorem~\ref{thm: inequalities with lvs and uvs}, we have 
\[\VPJ=\lvs=\uvs.\]
\end{definition}

 We now characterize equality in terms of a minimum-cardinality
set in $\mathcal D_c^-(\mathcal S)$ and a graph on the component indices. 

\begin{proposition}\label{prop: when S is v-combinatorial}
$\mS \in {\min}_{k}(G)$ is \textup{v}-{\em combinatorial} if and only if there exist $B \in \DCS$ of minimum cardinality and a simple graph $G'$ on $[k]$ such that:
\begin{enumerate}[label={\rm(\roman*)}]
\item\label{items:1} $G' \pitchfork \pms$ for every $\mS' \in \min(G)$ that is not nice with respect to $\mS$ and satisfies $\mS' \cap B = \emptyset$;
\item\label{items:2}
\[
|B \cap V_i| \geq \deg_{G'}(i) \quad \text{for all } i \in [k].
\]
\end{enumerate}
\end{proposition}

\begin{proof}

Suppose $\mathrm v^-_{\mathcal S}(G)
=\mathrm v^+_{\mathcal S}(G)$, and choose
$A\in\mathcal D_c^+(\mathcal S)$ such that
$|A_1|+2|A_2|=\mathrm v^+_{\mathcal S}(G)$.
Since $\operatorname{supp}(A)\in\mathcal D_c^-(\mathcal S)$,
\[
\mathrm v^-_{\mathcal S}(G)
\leq |\operatorname{supp}(A)|
\leq |A_1|+2|A_2|
=\mathrm v^+_{\mathcal S}(G)
=\mathrm v^-_{\mathcal S}(G).
\]
Thus equality holds throughout, and the elements of $A$
are pairwise disjoint. Set $B=\operatorname{supp}(A)$ and
$G'=G[A_2]$. The set $B$ has minimum cardinality.
Condition~\ref{items:1} follows from Definition~\ref{def: ECS}, and condition~\ref{item:2}
follows because distinct incident pairs use distinct vertices
of $B\cap V_i$.

Conversely, suppose $B,G'$ satisfy~\ref{items:1} and~\ref{items:2}.
For each $i$, assign distinct vertices of $B\cap V_i$ to
the edges of $G'$ incident with $i$; the degree condition
allows this assignment. Represent each edge $\{i,j\}$ by
the pair of assigned endpoints, and represent every unused
vertex of $B$ by a singleton. The resulting family $A$ has
support $B$, satisfies $G[A_2]=G'$, and has pairwise disjoint elements.
By~\ref{items:1} and Definition~\ref{def: ECS}, $A\in\mathcal D_c^+(\mathcal S)$.
Consequently,
\[
\mathrm v^+_{\mathcal S}(G)
\leq |A_1|+2|A_2|=|B|
=\mathrm v^-_{\mathcal S}(G).
\]
The reverse inequality follows from Theorem~\ref{thm: inequalities with lvs and uvs}.

\end{proof}

Proposition~\ref{prop: when S is v-combinatorial} gives an explicit combinatorial criterion for when the two
bounds coincide. We next derive a simple sufficient condition in terms of
the connected domination numbers of the components $G_i$.


\begin{theorem}\label{thm: all connected domination at least k-1}
If $\mS \in {\min}_{k}(G)$ and $\gamma_{c}(G_{i}) \geq k-1$ for all $i \in [k]$, then $\mS$ is $\mathrm{v}$-combinatorial and
\[
\VPJ = \lvs = \uvs.
\]
\end{theorem}

\begin{proof}
Let $B \subseteq [n] \setminus \mS$ be any element of $\DCS$ of minimal size, and let $\mathcal{K}_{k}$ denote the complete graph on $[k]$. We claim that $B$ and $\mathcal{K}_{k}$ satisfy conditions~\ref{items:1} and~\ref{items:2} of Proposition~\ref{prop: when S is v-combinatorial}, which implies that $\mS$ is $\mathrm{v}$-combinatorial.

To verify~\ref{items:1}, let $\mathcal S'$ be a non-nice cutset with
$\mathcal S'\cap B=\emptyset$. 
Then $\mathcal P_{\mathcal S'}$ is defined and is not the discrete partition $\{\{1\}, \ldots, \{k\}\}$. It follows that $\mathcal{K}_k$ has an edge inside one of its blocks,
 and hence $\mathcal{K}_{k} \pitchfork \pms$.

To verify~\ref{items:2}, note that since $B \in \DCS$, condition~(i) in Definition~\ref{def: DC-} implies that $B \cap V_{i} \in \mathcal{D}_{c}(G_{i})$ for all $i \in [k]$. Hence,
\[
|B \cap V_{i}| \geq \gamma_{c}(G_{i}) \geq k-1 = \deg_{\mathcal{K}_{k}}(i),
\]
where the second inequality follows from the assumption. This establishes~\ref{items:2}, completing the proof.
\end{proof}



Theorem~\ref{thm: all connected domination at least k-1} provides the local $\mathrm{v}$-numbers in a broad range of cases. We
next consider the cases $k=1,2,3$. Recall that we are assuming that~\eqref{hyp:H} holds.

\begin{theorem}\label{thm: k=1}
$\emptyset$ is $\mathrm{v}$-combinatorial. Moreover,
\[
\textup{v}_{P_{\emptyset}}(J_G) = \gamma_{c}(G).
\]
\end{theorem}

\begin{proof}
Applying Theorem~\ref{thm: all connected domination at least k-1} with $k=1$ shows that $\emptyset$ is $\mathrm{v}$-combinatorial. Consequently,
\[
\textup{v}_{P_{\emptyset}}(J_G) = \textup{v}_{\emptyset}^{-}(G)
= \min \{\, |A| : A \in \mathcal{D}_{c}^{-}(\emptyset) \,\}
\]
by Theorem~\ref{thm: explicit formula of lvs}. For $\mathcal{S} = \emptyset$, the definition of $\mathcal{D}_{c}^{-}(\emptyset)$ (see Definition~\ref{def: DC-}) reduces to $\mathcal{D}_{c}(G)$. Hence,
\[
\textup{v}_{P_{\emptyset}}(J_G)
= \min \{\, |A| : A \in \mathcal{D}_{c}(G) \,\}
= \gamma_{c}(G),
\]
which completes the proof.
\end{proof}

This result was previously established in \textup{\cite[Theorem~3.6]{ambhore2024v}}, \textup{\cite[Theorem~3.2]{jaramillo2024connected}}, and \textup{\cite[Theorem~3.9]{liwski2026v}}, each via different methods. The argument presented here provides a new conceptual proof, highlighting that $\emptyset$ is $\mathrm{v}$-combinatorial.

\begin{theorem}\label{thm: k=2}
Let $\mathcal S\in{\min}_2(G)$, with component
vertex sets $V_1,V_2$. Then $\mS$ is $\mathrm{v}$-combinatorial. Moreover,
\[
\VPJ = \gamma_{c}(G,\mS)
:=
\min \left\{ |A| : A \in \mathcal{D}_{c}(G,\mS) \right\},
\]
where
\[
\mathcal{D}_{c}(G,\mS)
=
\left\{
A \subseteq [n]\setminus \mS
\;\middle|\;
A \cap V_{i} \in \mathcal{D}_{c}\big(G_{V_{i}\cup \mS}\big)
\text{ for } i=1,2
\right\}.
\]
\end{theorem}

\begin{proof}
Since $\mS$ satisfies~\eqref{hyp:H}, we have $\gamma_{c}(G_{i}) \geq 1$ for $i=1,2$. Hence, applying Theorem~\ref{thm: all connected domination at least k-1} with $k=2$ shows that $\mS$ is $\mathrm{v}$-combinatorial. Consequently,
\[
\textup{v}_{P_{\mathcal{S}}}(J_G)
=
\textup{v}_{\mathcal{S}}^{-}(G)
=
\min \{\, |A| : A \in \DCS \,\}.
\]

By Definition~\ref{def: DC-}, a set $A \subseteq [n]\setminus \mS$ belongs to $\DCS$ if and only if $A \cap V_{i} \in \mathcal{D}_{c}(G_{i})$ for $i=1,2$, and every $s \in \mS$ is adjacent to at least one vertex in each $A \cap V_{i}$. This is equivalent to requiring that
\[
A \cap V_{i} \in \mathcal{D}_{c}\big(G_{V_{i}\cup \mS}\big)
\quad \text{for } i=1,2,
\]
and the result follows.
\end{proof}

This recovers \textup{\cite[Theorem~4.15]{liwski2026v}} under~\eqref{hyp:H}. We say that a graph is a {\em cone} if it contains a vertex adjacent to all other vertices.

\begin{theorem}\label{thm:k=3}
If $\mS \in {\min}_{3}(G)$ and none of the graphs $G_{i}$, for $i \in \{1,2,3\}$, is a cone, then $\mS$ is $\mathrm{v}$-combinatorial and
\[
\VPJ = \lvs = \uvs.
\]
\end{theorem}

\begin{proof}
Since none of the graphs $G_{i}$ is a cone, we have $\gamma_{c}(G_{i}) \geq 2$ for all $i \in [3]$. The claim then follows by applying Theorem~\ref{thm: all connected domination at least k-1} with $k=3$.
\end{proof}

\section{\texorpdfstring{A bounded gap between $\lvs$ and $\uvs$}
{Bounded gap}}\label{sec: gap between bounds}

In this section, we refine the bounds on the local $\mathrm{v}$-number by controlling the gap between $\lvs$ and $\uvs$. More precisely, we establish an inequality of the form
\[
\uvs \leq \lvs + C,
\]
for a suitable constant $C$, thereby showing that the local $\mathrm{v}$-number lies in a bounded interval, generalizing Theorem~\ref{thm: all connected domination at least k-1}. For the following result, we retain hypothesis~\eqref{hyp:H}. 

\begin{theorem}\label{thm: gap} Let $\mathcal S\in {\min}_{k}(G)$ be a cutset satisfying~\eqref{hyp:H}. We have \begin{equation}\label{gap}\uvs \leq \lvs+\sum_{i=1}^{k} \max\{k-1-\gamma_{c}(G_{i}),0\}, \end{equation} Consequently \[\VPJ\in [\lvs, \lvs+\sum_{i=1}^{k} \max\{k-1-\gamma_{c}(G_{i}),0\}].\] 
In particular, for $k\geq 2$, $\VPJ \in [\, \lvs,\; \lvs + (k-2)k]$.
\end{theorem}

\begin{proof}
Choose $B\in\mathcal D_c^-(\mathcal S)$ with
$|B|=\mathrm v^-_{\mathcal S}(G)$, and put
$b_i=|B\cap V_i|$. Each $b_i$ is positive and
$b_i\geq\gamma_c(G_i)$, since $B\cap V_i$ is connected and dominating in $G_i$.

For each edge $\{i,j\}$ of $\mathcal{K}_k$, choose one endpoint in
$B\cap V_i$ and one in $B\cap V_j$. Make the choices so that
the $k-1$ incidences at $i$ use exactly
$\min\{b_i,k-1\}$ distinct vertices: use distinct vertices
when possible, and repeat endpoints when necessary.
Let $C$ be the resulting family of pairs. Then
\[
G[C]=\mathcal{K}_k,\qquad |C|=\binom{k}{2},\qquad
|\operatorname{supp}(C)\cap V_i|=\min\{b_i,k-1\}.
\]
Define $A=A_1\amalg A_2$ by
\[
A_2:=C,
\qquad
A_1:=\bigl\{\{v\}:v\in
B\setminus\operatorname{supp}(C)\bigr\}.
\]
Then $\operatorname{supp}(A)=B$.
For every non-nice cutset $\mathcal S'$ avoiding $B$,
$\mathcal P_{\mathcal S'}$ is defined and is not discrete.
Therefore $\mathcal{K}_k$ has an edge inside one of its blocks, and
Definition~\ref{def: ECS} gives $A\in\mathcal D_c^+(\mathcal S)$.
Consequently,
\begin{align*}
\mathrm v^+_{\mathcal S}(G)
&\leq |A_1| + 2|A_2| \leq |B|-\sum_{i=1}^k\min\{b_i,k-1\}+k(k-1)\\
&=\mathrm v^-_{\mathcal S}(G)
  +\sum_{i=1}^k\max\{k-1-b_i,0\}\\
&\leq\mathrm v^-_{\mathcal S}(G)
  +\sum_{i=1}^k\max\{k-1-\gamma_c(G_i),0\}.
\end{align*}
The interval follows from Theorem~\ref{thm: inequalities with lvs and uvs}. For $k\geq2$,
$\gamma_c(G_i)\geq1$ gives the final bound $k(k-2)$.

\end{proof}

\section{Applications to cycles and trees}
\label{sec: examples}

In this section, we apply our framework to compute certain local $\mathrm{v}$-numbers for cycles and trees.

\subsection{Cycles}
\label{subsec: cycles}

Let $\mathcal S$ be a cutset of $C_n$ with
$|\mathcal S|\geq3$. We compute the two combinatorial bounds $\textup{v}_{\mathcal{S}}^{-}(C_n)$ and $\textup{v}_{\mathcal{S}}^{+}(C_n)$
directly from their transversal definitions, without assuming~\eqref{hyp:H}.
The empty cutset and cutsets that disconnect the graph into two connected components are treated
in \cite{ambhore2024v,jaramillo2024connected} and \cite{liwski2026v}, respectively.


\begin{notation}\label{notation c1 and c2}
Let $\mS$ be as above. We use the following notation.
\begin{itemize}
    \item The complement $[n]\setminus\mS$ is a disjoint union of intervals
    $I_1,\ldots,I_k$, where $I_j=[a_j,b_j]$ and
    $a_{j+1}=b_j+2$ for each $j\in[k]$. Note that
    $k=\size{\mS}$ and
    \[
        \{b_1+1,\ldots,b_k+1\}
        =
        \{a_1-1,\ldots,a_k-1\}
        =
        \mS.
    \]
    
    \item We set
    \[
        C_1(\mS)
        =
        \{j\in[k]:|I_j|=1\}.
    \]
\end{itemize}
The intervals and their endpoints are indexed cyclically:
$I_{k+1}=I_1$, and vertex labels are read modulo $n$.
\end{notation}

By Example~\ref{ex: Cn 1}, the vertices of $\mathcal S$ are pairwise
nonadjacent in $C_n$. Its complement is the union of the
intervals $I_1,\ldots,I_k$, where $k=|\mathcal S|$. By Example~\ref{ex: Cn 2}, we have
\[
\mathcal{D}(M(\mathcal{S}))
=
\{B\subseteq[n]:B\cap\mathcal{S}\neq\emptyset\}
\cup
\bigcup_{j=1}^{k}
\{B\subseteq[n]:|B\cap I_j|\geq2\}
\cup
\{B\subseteq[n]:|B|\geq3\}.
\]
By definition, $\TS$ is the collection of transversals of the hypergraph
\[
\left\{
\mathcal{D}(M(\mathcal{S}'))\setminus
\mathcal{D}(M(\mathcal{S}))
:
\mathcal{S}'\text{ is a cutset of }C_n,\ 
\mathcal{S}'\neq\mathcal{S}
\right\}.
\]

By \textup{\cite[Lemma~5.5]{liwski2026v}}, a subset
\[
A\subseteq
\left(\binom{[n]}{1}\cup\binom{[n]}{2}\right)
\setminus\mathcal{D}(M(\mathcal{S}))
\]
belongs to $\mathcal{T}(\mathcal{S})$ if and only if $\supp(A)=[n]\setminus\mathcal{S}$ and, for every $j\in[k]$, there are vertices
    $u\in I_j$ and $w\in I_{j+1}$ with $\{u,w\}\in A$.

Consequently, the upper combinatorial invariant is given by
\begin{equation}\label{eq: vp+ for Cn}
\textup{v}_{\mathcal{S}}^{+}(C_n)
=
\min\left\{
|A_1|+2|A_2|:A\in\mathcal{T}(\mathcal{S})
\right\}
=
n-\size{\mathcal{S}}+\size{C_1(\mathcal{S})}.
\end{equation}

We now focus on the lower invariant. 

\begin{lemma}\label{lem: vp- for Cn}
For a cutset $\mS$ of $C_n$ with $\size{\mS}\geq 3$, we have
\[
\pTS=\{[n]\setminus\mS\}.
\]
In particular,
\[
\textup{v}_{\mathcal{S}}^{-}(C_n)=n-\size{\mS}.
\]
\end{lemma}

\begin{proof}
By definition, $\pTS$ is the collection of transversals of the hypergraph
\[
\Delta_{\mathcal{S}}'
:=
\left\{
\mathcal{S}'\setminus\mathcal{S}
:
\mathcal{S}'\in\min(C_n)\setminus\{\mathcal{S}\},
\ \text{$\mathcal{S}'$ is nice with respect to $\mathcal{S}$}
\right\}.
\]
Let $A\subseteq[n]\setminus\mS$ be an arbitrary element of $\pTS$.
We show that $A=[n]\setminus\mS$ by proving that every
$i\in[n]\setminus\mS$ belongs to $A$. Using Notation~\ref{notation c1 and c2}, let
$i\in I_j=[a_j,b_j]$ for some $j\in[k]$. Suppose first that
$i\neq a_j,b_j$. 
Then $\mS':=\mS\cup\{i\}$ is a cutset of $C_n$ that is nice with respect to $\mS$. 
Since
$A\in\pTS$, we must have $A\cap(\mS'\setminus\mS)\neq\emptyset$.
As $\mS'\setminus\mS=\{i\}$, it follows that $i\in A$.

Suppose next that $i\in\{a_j,b_j\}$. Without loss of generality,
let $i=a_j$. If $a_j\neq b_j$, define $\mS':=(\mS\setminus\{i-1\})\cup\{i\}$.
Then $\mS'$ is a cutset of $C_n$ that is nice with respect to $\mS$,
and $\mS'\setminus\mS=\{i\}$.
Since $A\in\pTS$, we again obtain $i\in A$.

It remains to consider the case $a_j=b_j=i$. In this case, set $\mS':=(\mS\setminus\{i-1,i+1\})\cup\{i\}$.
The new vertex $i$ has only the two neighbors $i-1,i+1$,
both of which have been removed from the cutset. The remaining
vertices inherited from $\mathcal S$ are pairwise nonadjacent.
Also $|\mathcal S'|=|\mathcal S|-1\geq2$.
Thus $\mathcal S'$ is a cutset. The singleton interval disappears,
but distinct surviving original intervals are not joined, so
$\mathcal S'$ is nice with respect to $\mathcal S$. Since $\mS'\setminus\mS=\{i\}$, this implies $A\cap \{i\}\neq \emptyset$, and hence $i\in A$.


Therefore every element of $[n]\setminus\mS$ belongs to $A$, and hence
$A=[n]\setminus\mS$. Since $[n]\setminus\mS$ is itself a transversal,
we conclude that $\pTS=\{[n]\setminus\mS\}$.
The formula for $\textup{v}_{\mathcal{S}}^{-}(C_n)$ follows immediately.
\end{proof}

We can now combine the two bounds.

\begin{proposition}\label{prop: cycles}
Let $\mS$ be a cutset of $C_n$ with $\size{\mS}\geq3$. Then
\[
n-\size{\mS}
\leq
\textup{v}_{P_{\mathcal{S}}}(J_{C_n})
\leq
n-\size{\mS}+\size{C_1(\mS)}.
\]
In particular, if $C_1(\mS)=\emptyset$, then
\[
\textup{v}_{P_{\mathcal{S}}}(J_{C_n})=n-\size{\mS}.
\]
\end{proposition}

\begin{proof}
The result follows immediately from
Theorem~\ref{thm: inequalities with lvs and uvs}, together with
the formulas for $\textup{v}_{\mathcal{S}}^{+}(C_n)$ and
$\textup{v}_{\mathcal{S}}^{-}(C_n)$ established in
Equation~\eqref{eq: vp+ for Cn} and Lemma~\ref{lem: vp- for Cn},
respectively.
\end{proof}

These computations illustrate our lower-bound method for cutsets of cycles, including cases outside~\eqref{hyp:H}. The more precise cycle
results of \textup{\cite[Section~5]{liwski2026v}} use additional algebraic arguments.

\subsection{Trees}
\label{subsec: trees}

In this subsection, we apply the preceding results to determine
certain local $\mathrm v$-numbers of binomial edge ideals of trees.

Let $T$ be a tree on $[n]$, and let $\mathcal S$ be a cutset of $T$.
Recall that a subset of $V(T)$ is a cutset if and only if each
of its vertices has at least two neighbors outside the subset.
Throughout this subsection, we assume that $\mathcal S$
satisfies hypothesis~\eqref{hyp:H}. Write
\[
T_{[n]\setminus\mathcal S}=T_1\sqcup\cdots\sqcup T_k.
\]
Each $T_i$ is a noncomplete tree and therefore has at least
three vertices.

For a tree $H$, put
\[
\operatorname{iv}(H)
:=\{v\in V(H):\deg_H(v)\geq2\}.
\]
The internal vertices of each $T_i$ form its unique
inclusion-minimal connected dominating set. Consequently,
\[
\gamma_c(T_i)=|\operatorname{iv}(T_i)|,
\qquad
A_i\in\mathcal D_c(T_i)
\Longleftrightarrow
\operatorname{iv}(T_i)\subseteq A_i\subseteq V(T_i).
\]Combining Theorems~\ref{thm: inequalities with lvs and uvs} and~\ref{thm: gap} with Corollary~\ref{cor: sum of gamma c}, we obtain
\[
\begin{aligned}
\sum_{i=1}^{k}\size{\iv(T_i)}
&\leq \lvt
\leq \VPT\\
&\leq \uvt
\leq \lvt+\sum_{i=1}^{k}
\max\{k-1-\size{\iv(T_i)},0\}.
\end{aligned}
\]

In particular, if $\size{\iv(T_i)}\geq k-1$ for all $i\in[k]$,
then 
\[
\VPT=\lvt=\uvt.
\]

Moreover, by Theorem~\ref{thm: explicit formula of lvs}, $\lvt$ is equal to the minimum cardinality of a subset $A\subseteq[n]\setminus\mS$ satisfying the following conditions:
\begin{enumerate}[label={\rm(\roman*)}]
\item $A\cap V(T_i)\supseteq\iv(T_i)$ for every $i\in[k]$;

\item for every $s\in\mS$, there exist distinct $i,j\in[k]$ such that
\[
N_T(s)\cap A\cap V(T_i)\neq\emptyset
\quad\text{and}\quad
N_T(s)\cap A\cap V(T_j)\neq\emptyset.
\]
\end{enumerate}

We put together these observations into the following proposition.

\begin{proposition}\label{prop: trees 1}
Let $T$ be a tree and $\mS$ be a cutset satisfying~\eqref{hyp:H}, with components $T_1,\ldots,T_k$.
Suppose that $\size{\iv(T_i)}\geq k-1$ for all $i\in[k]$.
Then \[\VPT=\lvt=\uvt.\]
Furthermore, if every $s\in\mS$ is adjacent to vertices belonging to at least two of the sets $\iv(T_1),\ldots,\iv(T_k)$, then
\[
\VPT=\lvt=\uvt
=\sum_{i=1}^{k}\size{\iv(T_i)}.
\]
\end{proposition}

\begin{proof}
For a noncomplete tree component,
$\gamma_c(T_i)=|\operatorname{iv}(T_i)|$.
The first assertion follows from Theorem~\ref{thm: all connected domination at least k-1}.
For the second, put $B=\cup_i\operatorname{iv}(T_i)$.
The additional neighborhood condition makes
$B\in\mathcal D_c^-(\mathcal S)$.
Every member of this family contains each
$\operatorname{iv}(T_i)$, so $B$ has minimum cardinality.
Theorem~\ref{thm: explicit formula of lvs} and the first assertion give the formula.
\end{proof}

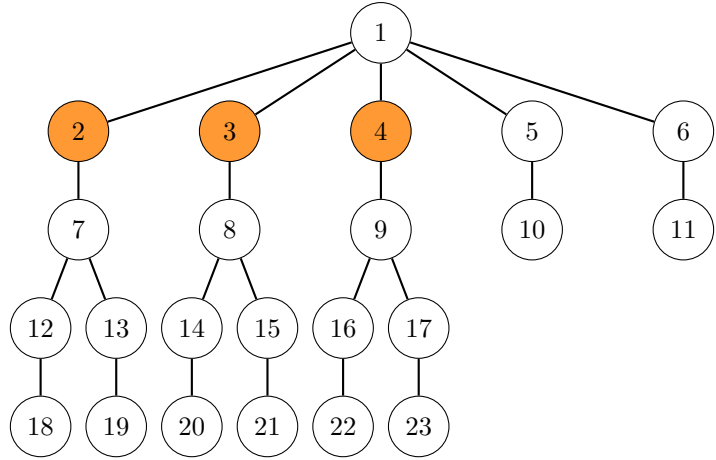
\begin{figure}[H]
    \centering
    \begin{tikzpicture}[
        every node/.style={
            circle,
            draw,
            minimum size=8mm,
            inner sep=0pt,
            font=\small
        },
        cutset/.style={
            fill=orange!80
        },
        edge/.style={
            draw,
            thick
        }
    ]

    \node (1) at (0,5) {1};

    \node[cutset] (2) at (-4,3.7) {2};
    \node[cutset] (3) at (-2,3.7) {3};
    \node[cutset] (4) at (0,3.7) {4};
    \node (5) at (2,3.7) {5};
    \node (6) at (4,3.7) {6};

    \node (7)  at (-4,2.4) {7};
    \node (8)  at (-2,2.4) {8};
    \node (9)  at (0,2.4) {9};
    \node (10) at (2,2.4) {10};
    \node (11) at (4,2.4) {11};

    \node (12) at (-4.5,1.1) {12};
    \node (13) at (-3.5,1.1) {13};
    \node (14) at (-2.5,1.1) {14};
    \node (15) at (-1.5,1.1) {15};
    \node (16) at (-0.5,1.1) {16};
    \node (17) at (0.5,1.1) {17};

    \node (18) at (-4.5,-0.2) {18};
    \node (19) at (-3.5,-0.2) {19};
    \node (20) at (-2.5,-0.2) {20};
    \node (21) at (-1.5,-0.2) {21};
    \node (22) at (-0.5,-0.2) {22};
    \node (23) at (0.5,-0.2) {23};

    \draw[edge] (1) -- (2);
    \draw[edge] (1) -- (3);
    \draw[edge] (1) -- (4);
    \draw[edge] (1) -- (5);
    \draw[edge] (1) -- (6);

    \draw[edge] (2) -- (7);
    \draw[edge] (3) -- (8);
    \draw[edge] (4) -- (9);
    \draw[edge] (5) -- (10);
    \draw[edge] (6) -- (11);

    \draw[edge] (7) -- (12);
    \draw[edge] (7) -- (13);
    \draw[edge] (8) -- (14);
    \draw[edge] (8) -- (15);
    \draw[edge] (9) -- (16);
    \draw[edge] (9) -- (17);

    \draw[edge] (12) -- (18);
    \draw[edge] (13) -- (19);
    \draw[edge] (14) -- (20);
    \draw[edge] (15) -- (21);
    \draw[edge] (16) -- (22);
    \draw[edge] (17) -- (23);

    \end{tikzpicture}
    \caption{A tree on $[23]$ with cutset $\{2,3,4\}$.}
    \label{fig:tree-23}
\end{figure}

\begin{example}\label{ex tree}
Let $T$ be the tree on $[23]$ shown in Figure~\ref{fig:tree-23}, and let $\mS:=\{2,3,4\}$.
The graph $T_{[23]\setminus\{2,3,4\}}$ has connected components $T_1,T_2,T_3,T_4$ with vertex sets
\[
\{7,12,13,18,19\},
\{8,14,15,20,21\},
\{9,16,17,22,23\},
\{1,5,6,10,11\},
\]
and sets of internal vertices
\[
\{\iv(T_1),\iv(T_2),\iv(T_3),\iv(T_4)\}
=
\big\{
\{7,12,13\},
\{8,14,15\},
\{9,16,17\},
\{1,5,6\}
\big\}.
\]
Thus, $\size{\iv(T_i)}=3=k-1$
for every $i\in[4]$.
Moreover, each $s\in\{2,3,4\}$ is adjacent to vertices belonging to at least two of the sets $\iv(T_1),\ldots,\iv(T_4)$. Hence, Proposition~\ref{prop: trees 1} gives
\[
\VPT=\sum_{i=1}^{4}\size{\iv(T_i)}=12.
\]
Since $|\operatorname{iv}(T)|=15$, this also gives
$\mathrm v(J_T)\leq12<15=\gamma_c(T)$.
Thus the general upper bound $\mathrm v(J_T)\leq\gamma_c(T)$
from \textup{\cite{jaramillo2024connected}} need not be sharp for trees.
\end{example}



We next consider singleton cutsets. Every internal vertex $i$
of a tree gives a cutset $\{i\}$. The following condition
determines its local v-number.

\begin{proposition}\label{prop: cutset of size 1}
Let $i\in[n]$ be such that every vertex in $N_T(i)\cup\{i\}$ has degree at least three. Then $\{i\}$ is $\mathrm{v}$-combinatorial, and
\[
\mathrm v_{P_{\{i\}}}(J_T)=\size{\iv(T)}-1.
\]
\end{proposition}

Each component of $T\setminus\{i\}$ contains a neighbor of
$i$ with at least two further neighbors. It therefore has
at least three vertices and satisfies~\eqref{hyp:H}. To prove this proposition, we first determine the lower bound $\mathrm{v}_{\{i\}}^{-}(T)$.

\begin{lemma}\label{lem: cutset of size 1}
Under the hypotheses of Proposition~\ref{prop: cutset of size 1}, the set
$\operatorname{iv}(T)\setminus\{i\}$ is the unique
inclusion-minimal element of $\mathcal T'(\{i\})$. 
 Consequently,
\[
\mathrm{v}_{\{i\}}^{-}(T)=\size{\iv(T)}-1.
\]
\end{lemma}

\begin{proof}
We first show that $\iv(T)\setminus\{i\}\in\pti$. Let $\mS\neq\{i\}$ be a cutset that is nice with respect to $\{i\}$. Suppose, for contradiction, that $\mS\cap\bigl(\iv(T)\setminus\{i\}\bigr)=\emptyset$.
Since every vertex in $[n]\setminus\iv(T)$ has degree one, no such vertex can belong to a cutset. Hence $\mS\subseteq\{i\}$.
As $\mS\neq\{i\}$, this forces $\mS=\emptyset$. However, $\emptyset$ is not nice with respect to $\{i\}$, a contradiction. Therefore, $\iv(T)\setminus\{i\}\in\pti$.

It remains to show that every $A\in\pti$ contains $\iv(T)\setminus\{i\}$. Let $j\in\iv(T)\setminus\{i\}$. 
We claim that $\{i,j\}$ is a cutset of $T$.
If $j\notin N_T(i)$, then both $i$ and $j$ have at
least two neighbors outside $\{i,j\}$.
If $j\in N_T(i)$, then both vertices have degree at
least three by hypothesis, so again each has at least
two neighbors outside $\{i,j\}$.
The characterization of cutsets of trees recalled
at the beginning of this subsection therefore shows
that $\{i,j\}$ is a cutset.
Moreover, it is nice with respect to $\{i\}$.
Since $A\in\pti$, we must therefore have $A\cap\{i,j\}\neq\emptyset$.
By definition, $A\subseteq[n]\setminus\{i\}$, and hence $j\in A$. Since $j\in\iv(T)\setminus\{i\}$ was arbitrary, we conclude that $\iv(T)\setminus\{i\}\subseteq A$.
Thus, $\iv(T)\setminus\{i\}$ is the unique inclusion-minimal element of $\mathcal{T}'(\{i\})$, completing the proof.
\end{proof}

\begin{proof}[Proof of Proposition~\ref{prop: cutset of size 1}]
By Lemma~\ref{lem: cutset of size 1}, $\mathrm{v}_{\{i\}}^{-}(T)=\size{\iv(T)}-1$.
Thus, by Theorem~\ref{thm: inequalities with lvs and uvs}, it suffices to show that $\{i\}$ is $\mathrm{v}$-combinatorial. We apply Proposition~\ref{prop: when S is v-combinatorial}. We need to exhibit a set $B\in\pti$ of minimal size and a graph $G'$ on $[k]$ such that:
\begin{enumerate}[label={\rm(\roman*)}]
\item $G'\pitchfork\pms$ for every $\mS'\in\min(T)$ that is not nice with respect to $\{i\}$ and satisfies $\mS'\cap B=\emptyset$;
\item
\[
|B\cap V(T_j)|\geq\deg_{G'}(j)
\qquad\text{for all }j\in[k].
\]
\end{enumerate}

Set $B:=\iv(T)\setminus\{i\}$.
By Lemma~\ref{lem: cutset of size 1}, this set has minimal size in $\pti$. 
Since leaves belong to no cutset of a tree, every cutset
avoiding $B=\operatorname{iv}(T)\setminus\{i\}$ is contained
in $\{i\}$.
Thus the only such cutsets are $\emptyset$ and $\{i\}$.
The latter is nice with respect to itself, whereas
$\emptyset$ is not nice with respect to $\{i\}$.
Therefore $\emptyset$ is the only non-nice cutset
avoiding $B$ that must be considered in condition~(i). 
Hence, if $G'$ is the graph on $[k]$ with edge set $E(G')=\{\{1,2\}\}$,
then $G'\pitchfork\mathcal{P}_{\emptyset}$.
It remains to verify condition~{\rm(ii)}. Since every neighbor of $i$ has degree at least three in $T$, each such neighbor is an internal vertex of $T\setminus\{i\}$. Consequently, every connected component $T_1,\ldots,T_k$ of $T\setminus\{i\}$ contains at least one internal vertex. Thus,
\[
|B\cap V(T_j)|\geq1
\qquad\text{for all }j\in[k].
\]
Since the graph $G'$ has maximum degree one, we have
\[
|B\cap V(T_j)|\geq\deg_{G'}(j)
\qquad\text{for all }j\in[k].
\]
Therefore, both conditions of Proposition~\ref{prop: when S is v-combinatorial} are satisfied, and $\{i\}$ is $\mathrm{v}$-combinatorial.

Hence, by Theorem~\ref{thm: inequalities with lvs and uvs},
\[
\mathrm v_{P_{\{i\}}}(J_T)
=
\mathrm{v}_{\{i\}}^{-}(T)
=
\size{\iv(T)}-1.
\]
\end{proof}

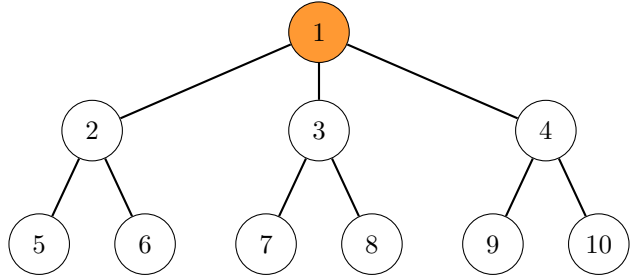
\begin{figure}[H]
    \centering
    \begin{tikzpicture}[
        every node/.style={
            circle,
            draw,
            minimum size=8mm,
            inner sep=0pt,
            font=\small
        },
        cutset/.style={
            fill=orange!80
        },
        edge/.style={
            draw,
            thick
        }
    ]

    \node[cutset] (1) at (0,4) {1};

    \node (2) at (-3,2.7) {2};
    \node (3) at (0,2.7) {3};
    \node (4) at (3,2.7) {4};

    \node (5)  at (-3.7,1.2) {5};
    \node (6)  at (-2.3,1.2) {6};
    \node (7)  at (-0.7,1.2) {7};
    \node (8)  at (0.7,1.2) {8};
    \node (9)  at (2.3,1.2) {9};
    \node (10) at (3.7,1.2) {10};

    \draw[edge] (1) -- (2);
    \draw[edge] (1) -- (3);
    \draw[edge] (1) -- (4);

    \draw[edge] (2) -- (5);
    \draw[edge] (2) -- (6);
    \draw[edge] (3) -- (7);
    \draw[edge] (3) -- (8);
    \draw[edge] (4) -- (9);
    \draw[edge] (4) -- (10);

    \end{tikzpicture}
    \caption{A tree on $[10]$ with cutset $\{1\}$.}
    \label{fig:tree-10}
\end{figure}

\begin{example}
Let $T$ be the tree on $[10]$ shown in Figure~\ref{fig:tree-10}, and let $\mS=\{1\}$. Since the vertex $1$ and all of its neighbors have degree at least three, Proposition~\ref{prop: cutset of size 1} implies that $\mS$ is $\mathrm{v}$-combinatorial and
\[
\mathrm v_{P_{\{1\}}}(J_T)=\size{\iv(T)}-1=3.
\]
\end{example}

The preceding examples illustrate how a suitable choice of cutset can improve the global upper bound. In the final section we discuss a possible connection with regularity.

\section{Comparison with regularity and open questions}\label{regularity}

It was conjectured in \textup{\cite[Conjecture~5.3]{ambhore2024v}} that $\gamma_{c}(G)\leq \reg(S/J_G)$ for every connected graph $G$, which would imply $\mathrm{v}(J_G)\leq \reg(S/J_G)$. 
Known cases include trees \textup{\cite[Corollary~3.7]{jaramillo2024connected}}, and, more generally, chordal graphs, as well as whisker graphs, \textup{\cite[Section~4]{ambhore2024v}}. An analogous comparison arises for Lov\'asz--Saks--Schrijver ideals. Over an algebraically closed field, the inequality 
\[\vn(L_G^{\mathbb{K}}(d)) \leq \reg(R/L_G^{\mathbb{K}}(d))\]
was established for forests and every $d\geq 1$ in \textup{\cite[Theorem~4.14]{kumar2026mathrm}}, where $R$ is the polynomial ring containing $L_G^{\mathbb{K}}(d)$.

Our upper bounds suggest allowing the minimal prime to vary. Set
\[U(G):=\min_{\mathcal{S}\in \min(G)}\uvs.\]
By Theorem~\ref{thm: inequalities with lvs and uvs} and~\ref{thm: k=1}, for every connected noncomplete graph 
\[\vn(J_{G})\leq U(G)\leq \gamma_c(G).\]
The second inequality can be strict, as Example~\ref{ex tree} shows. Thus $U(G)$ may provide a sharper bound when comparing the global $\vn$-number with the regularity. These observations suggest the following conjecture.

\begin{conjecture}
For every connected graph $G$,
\[U(G)\leq \reg(S/J_G).\]
\end{conjecture}

This conjecture is implied by $\gamma_c(G)\leq \reg(S/J_G)$ and would itself imply $\vn(J_{G})\leq \reg(S/J_G)$.

A possible approach is to compare $U(G)$ with the combinatorial invariant $\nu(G)$ introduced in \textup{\cite[Definition~29]{Laclair}}. It satisfies
\[\nu(G)\leq \reg(S/J_G),\]
with equality for closed graphs, block graphs, and bipartite graphs whose binomial edge ideals are Cohen--Macaulay; see \textup{\cite[Theorem~30, Corollaries~33 and~37, and Theorem~64]{Laclair}}.

In particular, for block graphs and closed graphs, which are chordal, the connected domination bound established in \textup{\cite[Theorem~4.5]{ambhore2024v}} together with \textup{\cite[Corollary~33 and Theorem~64]{Laclair}}, gives
\[U(G)\leq \gamma_c(G)\leq \reg(S/J_G)=\nu(G).\]
For these families the empty cutset therefore suffices. Allowing other cutsets leads to the following question.

\begin{question}\label{question nu}
For which connected graphs $G$ does the inequality $U(G)\leq \nu(G)$ hold? 
\end{question}

For every graph satisfying this inequality, one has 
\[\vn(J_{G})\leq U(G)\leq \nu(G)\leq \reg(S/J_{G}).\]


\printbibliography

\medskip
{\footnotesize\noindent {\bf Authors' addresses}
\medskip

\noindent{Emiliano Liwski, 
KU Leuven}\hfill {\tt  emiliano.liwski@kuleuven.be}

\end{document}